\documentclass[11pt]{amsart}

\usepackage{graphics}

\usepackage{amsmath,amsfonts,amssymb,mathrsfs}
\usepackage{graphics}
  \usepackage{epsfig}

\newtheorem{theorem}{Theorem}[section]

\newtheorem{remark}[theorem]{Remark}
\newtheorem{lemma}[theorem]{Lemma}
\newtheorem{definition}[theorem]{Definition}

\newcommand\ws{{\widetilde {\mathcal{L}}}}
\newcommand\mcs{\mathcal{S}}

\newcommand\wv{{\widetilde V}}
\newcommand\ep{{\varepsilon}}

\newcommand\ph{\varphi}
\newcommand\de{\delta}
\newcommand\zu{[0,1]}
\newcommand\zuo{[0,1]}

\newcommand\ms{\medskip}
\newcommand\sk{\smallskip}

\newcommand{\R}{\mathbb{R}}
\newcommand{\N}{\mathbb{N}}
\newcommand{\Z}{\mathbb{Z}}

\begin{document}
\title{A localized  Jarnik-Besicovitch theorem}

\author{Julien Barral \and St\'ephane Seuret }         

\address{Julien Barral, Projet SISYPHE, INRIA Rocquencourt, B.P. 105, 78153 Le Chesnay Cedex,
France  - St\'ephane Seuret,   LAMA, CNRS UMR 8050, Universit\'e Paris-Est,  
  61 avenue du G\'en\'eral de Gaulle, 
94 010 CR\'ETEIL Cedex, France}
\date{Received: date / Revised version: date}

\begin{abstract}
Fundamental questions in Diophantine approximation are related to the  Hausdorff dimension of sets of the form  $\{x\in \R: \de_x = \de\}$, where $\de \geq 1$ and $\de_x$ is the Diophantine approximation rate of an irrational number
$x $. We go beyond the classical results by computing the Hausdorff dimension of the sets   $\{x\in\R: \de_x =f(x)\}$, where $f$ is a continuous function. Our theorem   applies to the study of the approximation rates by various approximation families. It also applies to functions $f$ which are continuous outside a set of prescribed Hausdorff dimension. \end{abstract}

\maketitle

\section{Introduction}

Let $(x_n)_{n\ge 1}$ be a sequence    in a $\sigma$-compact metric space $(E,d)$, and  $(r_n)_{n\ge 1}$ be a non-increasing sequence  of real numbers converging to 0 when $n$ tends to infinity. Let $\ph:\R_+ \to \R_+$ be a positive non-decreasing continuous mapping with $\varphi(0)=0$. The   set
\begin{equation}
\label{defsdelta} 
  \mathcal{L}(\ph)   =  \{x\in E: \ d(x,x_n) \leq  \ph(r_n)   \mbox{ for infinitely many integers } n\}
\end{equation}
 contains the elements of $E$ that are infinitely often well-approximated, at rate $\ph$, by the points $x_n$ relatively to the radii $r_n$. This set can be rewritten as a limsup set:
$$\, \mathcal{L} (\ph) =\limsup_{n\to\infty} B\big (x_n,\varphi(r_n)\big ).$$  
($B(x,r)$ stands for the ball of centre $x$ and radius $r$.)
The values of the Hausdorff dimensions (or the Hausdorff measure associated with convenient gauge functions) of sets of the form $\, \mathcal{L}(\ph)$ provide us with  a fine description of the geometrical distribution in $E$ of the sequence   $(x_n)_{n\geq 1}$.  

Such limsup sets arise naturally in  Diophantine approximation theory in $\mathbb{R}^d$ (see \cite{khin,JARNIK,khin2,BESIC,guting,BUGEAUD,BUGEAUD2,BSUbiquity,BSUbiquity2,VEL4} among many references), and more generally in Diophantine approximation problems in limit sets of groups or  in Julia sets of rational maps \cite{DODVEL,STRAT,HILL,HILL2}. They also   appear in mathematical physics and dynamical systems when studying resonance problems \cite{arnold,Petro1,Petro2,Ptas,Beres3}, and    when measuring the distribution of H\"older singularities of measures and functions \cite{JAFFARD,JAFFARD1,JAFFARD2,Falc,BS1,BS2}.

\ms

We  denominate the sequence of couples  $\mathcal{S}= \big ((x_n,r_n)\big )_{n\ge 1}$ as an {\em approximation system} (or simply a {\em system}) in $E$. Standard Diophantine approximation deals with the approximation of real numbers by the  system $  \big((p/q, 1/q^2)\big)_{q\geq 1, \, p\in \Z }$.   

\ms

The  mappings  $\ph$    of the form $\ph_\de: r\mapsto r^{\delta}$ (for $\delta >0$) are particularly relevant. Hence, denoting $\, \mathcal{L}({\ph_\de})$   simply by $ \mathcal{L}_{\de}$, we consider
\begin{equation}
\label{defjarnik0} 
 \mathcal{L}_\de   =  \left\{x \in E: d(x,x_n) \leq  r_n^{\de} \ \mbox{ for infinitely many integers $n$}\right\}.
\end{equation}
The sets $( \mathcal{L}_\de)_{\de > 0}$ form a non-increasing family of sets.  This   property allows us to classify the elements $x$ of $E$ according to their   {\em approximation rate}   by the system $ \mathcal{S}= \big ((x_n,r_n)\big )_{n\ge 1}$. This approximation rate is  defined for $x\in E$ by  (we use the convention that $\sup \emptyset =0$)
 $$\de_x=\sup \{\de: x\in  \mathcal{L}_\de\},$$
and for $\de\geq 0$,  one is naturally interested in  the set $\widetilde { \mathcal{L}}_\de$  of points which have  approximation rate $\de$:
\begin{equation}
\label{defst}
\widetilde { \mathcal{L}}_\de = \{x \in E: \, \delta_x=\delta\}. \ \
\end{equation}
We emphasize the following embedment properties between the sets $ \mathcal{L}_\de$ and $\ws_\de$: For $\de>0$, 
\begin{equation}
\label{embed}
 \mathcal{L}_\de  =    \bigcup_{\delta'\geq \delta}\widetilde { \mathcal{L}}_{\de'}    \ \ \mbox { and } \  \  \ws_\de =  \bigcap_{\de'\geq \de} { \mathcal{L}}_{\de'} \setminus \bigcup_{\de'<\de}  \mathcal{L}_{\de'} .\end{equation}
The dimension problems related with the sets $\mathcal{L}_{\de}$ are also relevant for the sets  $\ws_\de$. Hence, given $\delta > 0$,  it is natural to question   the non-emptiness of $\ws_\de$,  and the value of the Hausdorff dimension of $\ws_\de$, and  the existence of  gauge functions $\zeta$ for which the corresponding Hausdorff measure $\mathcal{H}^\zeta(\ws_\de)$ is null, positive and finite, or infinite (
$\mathcal{H}^\xi $ stands for the generalized Hausdorff measure associated with the gauge function $\xi$, see Section \ref{Hausdorff}
). The first investigations on this subject have led to the celebrated Jarnik-Besicovitch theorem: if  the system $\mcs$ is the rational system $  \big((p/q, 1/q^2)\big)_{q\geq 1, \, p\in \Z }$,  then for every $\de \geq 1$, $\dim_{\mathcal{H}} \, \mathcal{L}_\de = 1/\de$ ($\dim_{\mathcal{H}}$ stands for the Hausdorff dimension). 
 
\ms

Note that in the case of the approximation by rational numbers, $ \mathcal{L}_1 = \R$, since by a famous Dirichlet's result, for every $x\in \R$, the event $x\in B(p/q, 1/q^2)$ occurs for infinitely many integers $q\geq 1$.
In fact, in many situations, the system $\mcs$  is chosen so that  the set
\begin{equation}
\label{eq000}
\mathcal{L}_1=\limsup_{n\to\infty}B(x_n,r_n) \ \ \mbox{is of full $m$-measure in $E$,}
\end{equation}
where $m$ is a probability measure on $E$ enjoying nice scaling properties:
\begin{itemize}
\item
  In \cite{JAFFLACUNARY,DODVEL}, $m$ is {\em uniformly distributed}: There exist $r_0>0$ and $C>1$ such that $C^{-1} r^{\dim_{\mathcal{H}}E}\le m(B(x,r))\le Cr^{\dim_{\mathcal{H}} E}$ for all $x\in E$ and $r\in (0,r_0]$. Observe that if such a measure $m$ exists on $E$, then  ${\dim_{\mathcal{H}}}\,  (m)=  {\dim_{\mathcal{H}}}\,  (E)$ (${\dim_{\mathcal{H}}}\,  (m)$ is the Hausdorff dimension of the measure $m$).
  \item
   In \cite{BSUbiquity}, $m$ is supposed to possess deterministic or statistical self-similarity properties, which imply the weaker property: There exists $\delta\in (0,\dim_{\mathcal{H}}\, E]$ such that $\lim_{r\to 0^+} \frac{\log(m(B(x,r))}{\log (r)}={\dim_{\mathcal{H}}}\,  (m)$ for $m$-almost every $x$. 
   \end{itemize}
   In these contexts, it can be  proved that for every  $\de\ge1$, ${\dim_{\mathcal{H}}}\,  \mathcal{L}_\de\ge {\dim_{\mathcal{H}} }\, (m)/\delta$. Moreover, \, there exists a gauge function \, $\xi:\R_+\to \R_+$ \, satisfying  $\lim_{r\to 0^+} \frac{\log \xi(r)}{\log r}   =\dim_{\mathcal{H}}\,  (m)/\de$ and $\mathcal{H}^\xi ( \mathcal{L}_\de)>0$.
Such theorems are referred to as {\em ubiquity} results, and the literature on ubiquity properties is numerous.

\ms 

From now on, we focus on $\R^d$, with $d\geq 1$, and more precisely, for obvious periodicity reason, on $E=\zuo^d$.  

As mentioned previously, the authors of \cite{JAFFLACUNARY,DODVEL} obtained  the following ubiquity theorem which treats the case where the measure $m$ in (\ref{eq000}) equals $\ell$,    the Lebesgue measure in $\R^d$. Theorem \ref{thubi} establishes  an extension of the famous Jarnik-Besicovitch theorem for Diophantine approximation by rational numbers.
\begin{theorem}
\label{thubi}
Let $\mcs=\Big((x_n,r_n)\Big)_{n\geq 1}$ be a system in $\zuo^d$.  If 
\begin{equation}
\label{cond1}
\ell \, \Big( \mathcal{L}_1 \Big)  = 1,   \end{equation}
then for every $\de \geq 1$, ${\dim_{\mathcal{H}}} \,   \mathcal{L}_\de \geq  d/\delta$.

In addition, \, there exists a gauge function \, $\xi:\R_+\to \R_+$ satisfying  \, $\lim_{r\to 0^+} \frac{\log \xi(r)}{\log r} =d/\de$ and $\mathcal{H}^\xi ( \mathcal{L}_\de)>0$.
\end{theorem}

The second part of Theorem \ref{thubi} is crucial, since it deals with Hausdorff  measures (and not only with the Hausdorff dimension). It makes it possible to replace  the set  $ \mathcal{L}_\de$   by $\ws_\de$ in the statement of Theorem \ref{thubi}, provided that  the balls $B(x_n,r_n)$ with comparable diameters  do not overlap excessively. This occurs  when there exists  an integer $N>0$ such for all $j\ge 0$, each element  $x\in \zuo^d$  belongs to at most $N$ balls $B(x_n,r_n)$ such that $2^{-j-1}\le r_n\le 2^{-j} $ (heuristically, the elements of $\zuo^d$ are covered "economically" by the balls $B(x_n,r_n)$). Such a property is a specific case of the {\it weak redundancy} property $\mathcal{C}_1$, which will be defined in Section \ref{WRe}.  In this case, we thus have for all $\de\ge 1$
\begin{equation}
\label{geom}
{\dim_{\mathcal{H}}} \, \mathcal{L}_\de =    {\dim_{\mathcal{H}}} \, \widetilde {\mathcal{L}}_\de= \frac{d}{\de}.
\end{equation}
This two-sided  equality  contains two results: the non-emptiness of $\ws_\de$ and the value of its Hausdorff dimension.  It (\ref{geom}) holds  in $\mathbb{R}$ when  considering the  ``rational" system  $\mathcal{R}=\big((p/q, 1/q^2)\big)_{q\geq 1, \, 0<p<q }$ or other systems of points (for instance obtained as Poisson point processes in the upper-half-plane, see Section \ref{secexa} for details and further examples).
 
\ms

Finally, observe that for a given system $\mathcal{S}= \big ((x_n,r_n)\big )_{n\ge 1}$, the covering property
\begin{equation}
\label{optim}
\limsup_{n\to\infty} B(x_n,r_n) =(0,1)^d
\end{equation}
implies (\ref{cond1}), and   the  corresponding approximation rates satisfy
$$\de_x\ge 1 \ \  \mbox{ for all $x \in (0,1)^d$.}
$$
Therefore, the associated sets $\ws_\de$ provide us with  a classification of all the elements of  $\zuo^d$ with respect to their approximation rates (those associated with $\mathcal{S}$).

\ms

 In this article we replace the (constant) approximation rate $\delta$ in (\ref{defst}) by $f(x)$, where $f$ is  a {\em continuous function}. We are thus looking for elements $x\in \zuo^d$ whose {\em approximation rate by some system $\mcs$ depend on $x$ via the function $f$}. Hence, Jarnik-Besicovitch's Theorem and the result on the Hausdorff dimension of $\widetilde{\mathcal{L}}_\de$ in Theorem \ref{thubi} will be viewed as Theorem \ref{main} in the special case where $f$ is a constant function.

Let us   state our main theorem. Conditions $\mathcal{C}_1$ and $\mathcal{C}_2$ will be explained later. 
\begin{theorem}
\label{main}
Consider the system  $\mathcal{S}= \big ((x_n,r_n)\big )_{n\ge 1} $, where $(x_n)_{n\ge 1}$  is a sequence of elements of  $(0,1)^d$  and  $(r_n)_{n\ge 1}$ is a non-increasing sequence of real numbers converging to 0 when $n$ tends to infinity.

Assume that (\ref{optim}), $\mathcal{C}_1$ and $\mathcal{C}_2$ hold.

Let  $\Omega$ be a non-empty compact subset of $(0,1)^d$, such that $\overline{\overset{\circ}{\Omega}} =\Omega$.

Let $f:(0,1)^d \to [1,+\infty)$ be a continuous function.

Consider the subsets of  $ \zuo^d $ defined by
\begin{eqnarray}
\, \mathcal{L}(\Omega,f)&  = &  \left\{x\in\Omega: \ \delta_x \geq f(x) \right\}\\
\widetilde{ \mathcal{L}}(\Omega,f) & = & \left\{x\in\Omega: \ \delta_x  =  f(x) \right\}
\end{eqnarray}

 The sets $\mathcal{L}(\Omega,f)$ and $\widetilde { \mathcal{L}}(\Omega,f) $ are dense in $\Omega$ and we have 
\begin{equation}
\label{result}
{\dim_{\mathcal{H}}} \, \mathcal{L}(\Omega,f)= {\dim_{\mathcal{H}}}\, \ws(\Omega,f)     = \frac{d}{\min \{f(x): x\in \Omega\}}.
\end{equation}
\end{theorem}

\noindent The function $f$ ranges over $[1,+\infty)$, since $\de_x$ is always larger than $1$. 

\ms

As stated above, formula (\ref{geom}) shall now be seen as a particular case of (\ref{result}). As in relation (\ref{geom}),    (\ref{result}) contains  several results: the non-emptiness of $\ws(\Omega,f)$, the equality between the Hausdorff dimensions of $\, \mathcal{L}(\Omega,f)$ and $\ws(\Omega,f)$, and the value of this dimension. 

 \ms
 
 The key point is that  conditions $\mathcal{C}_1$ and $\mathcal{C}_2$ hold for many classical systems   arising in ubiquity  and number theory.  In Section \ref{secexa}  we prove that Theorem \ref{main} applies to the  Diophantine approximation by dyadic numbers $\mathcal{D}= \Big(({\bf k}\cdot 2^{-j}, 2^{-j})\Big)_{j\geq 1, {\bf k}\in \{0,1,...,2^j-1\}^d}$, to the Diophantine approximation by  rational numbers $\mathcal{R} = \big ((p/q, 1/q^{2}) \big)_{q\geq 1, \, 0\le p\le q-1}$, to the so-called "inhomogeneous" Diophantine approximation by the system $\mathcal{I}= \Big(\big(\{n\alpha\}, \frac{1}{n}\big)\Big)_{n\geq 1}$, (where $\alpha$ is an irrational number whose approximation rate by the rational system $\mathcal{R}$ equals 2), and to the approximation rates by Poisson point processes $\mathcal{P}$.

%
%

\begin{center}
\begin{figure}
\begin{center}\label{fig1}
Figures are available on our web site.
\caption{Geometric representation of $T(\Omega,f)$ and $\widetilde T(\Omega,f)$ for the constant function $f(x)=\de$ on the left figure, for a typical continuous function $f$ on the right figure.}
\end{center}\end{figure}
\end{center}

Equality (\ref{result}) can be interpreted geometrically. 
Consider the subsets of  $ \R^d\times \R$ defined by
\begin{eqnarray*}
T(\Omega,f)&  = & \left\{(x,\delta_x): x\in\Omega \ \mbox{ and } \ \delta_x \geq f(x) \right\}\\
\widetilde T(\Omega,f) & = & \left\{(x,\delta_x): x\in\Omega \ \mbox{ and } \  \delta_x  = f(x) \right\}
\end{eqnarray*}
Then 
$\, \mathcal{L}(\Omega,f) $ and $\ws(\Omega,f) $ are respectively  the  natural projections of $T(\Omega,f)$  and $\widetilde T(\Omega,f)$ on $\R^d$. Theorem \ref{main} asserts that the ``frontier" of $T(\Omega,f)$, $\widetilde T(\Omega,f)$, is non empty and that the  projections of $T(\Omega,f)$ and $\widetilde T(\Omega,f)$ on $\R^d$ are both  dense in $\Omega$  and have   same Hausdorff dimension.
 
 \medskip
 
 Changing our standpoint, Theorem \ref{main} makes it possible to answer the following questions: 
 $$
 \mbox{Are there  real  numbers  $x\in \zu$   satisfying $\de_x =  1+x$? \ \  $\de_x=1/x$?}
 $$
This question is of course not reachable via Jarnik's result, for which the approximation rate is a fixed number $\de\geq 1$, independent of $x$. Moreover, it seems non-trivial (though possible) to explicitly construct an irrational number $x\in [0,1]$ such that $\delta_x=1+x$.  Theorem \ref{main} implies for instance that,  provided that a system  $\mathcal{S}$ in $\zu$ satisfies (\ref{optim}), $\mathcal{C}_1$ and $\mathcal{C}_2$, then   
\begin{itemize}
\item 
for every real numbers  $0< a <b < 1$,    
$$
   {\dim_{\mathcal{H}}} \ \Big\{x \in [a,b] : \, \de_x = 1+x\Big\}  =\frac{1}{1+a},$$

\item
for every real numbers  $0< a < b < 1$, for every $ \alpha \geq 1$,
$$ {\dim_{\mathcal{H}}} \,\Big\{x \in [a,b]: \de_x  = \frac{\alpha}{x} \Big\}    = \frac{b}{\alpha },$$
\item
and if  $ [a,b]\subset \displaystyle \Big [\frac{1}{6},\frac{5}{6}\Big ]$, then
 \begin{eqnarray*}
  {\dim_{\mathcal{H}}} \, \Big\{x \in [a,b]:\,  \de_x  = 2\sin({\pi x}) \Big\} \! & = \! & \displaystyle \frac{1}{\min \big (2\sin(\pi a),2\sin(\pi b)\big )}.
\end{eqnarray*}

\end{itemize}
 
In the  above equalities, the dimensions depend on the  range of $x$. This was expected, since the conditions we impose on $x$ depend on the non-constant continuous function $f$. 

\ms 
 
In order to prove Theorem \ref{thubi} and the equality between the Hausdorff dimensions of the (classical) sets ${\mathcal{L}}_\de$ and $\ws_\de$ defined respectively in  (\ref{defjarnik0}) and (\ref{defst}),   the usual method consists in  constructing iteratively a Borel probability measure $m_\delta$ of Hausdorff dimension larger than or equal to $d/\delta$ supported by ${\mathcal{L}}_\de$. Then,  recalling that  ${\dim_{\mathcal{H}}}\,  \ws_{\de'} < d/\de$ when $\de'>\de$,  we deduce from (\ref{embed})  that $m_\de(\ws_\de) =1$. Hence $m_\de$ is supported by $\ws_\de$ and ${\dim_{\mathcal{H}}} (\ws_\de) =d/\de$. Moreover, $m_\de$ can be chosen as the Hausdorff measure associated with a suitable  gauge function $g$  satisfying $\lim_{r\to 0^+}\log (g(r))/\log (r)=1/\delta$. 

As shall be explained soon, this approach is   inappropriate in the context of Theorem \ref{main}. First, observe that $\, \mathcal{L}(\Omega,f)$ and $\widetilde { \mathcal{L}}(\Omega,f) $  cannot be written as limsup sets. Nevertheless we still need to   construct probability measures with support contained in $\widetilde { \mathcal{L}}(\Omega,f)$. This set   is dense in $\Omega$ (like  $\widetilde {\mathcal{L}}_\delta$ in the introduction), but in general it is mostly localized around those elements of $\Omega$ at which $f$ reaches its minimum. This induces that in general, if $B$ is a non-trivial closed ball inside $\Omega$, then we have ${\dim_{\mathcal{H}}} \, ( \widetilde { \mathcal{L}}(\Omega,f)\cap B)={\dim_{\mathcal{H}}}\,  \widetilde{ \mathcal{L}} (\Omega,f)$ only if $f$ reaches its minimum over $B$. This constitutes a notable difference with the sets $\ws_\delta$, for which ${\dim_{\mathcal{H}}}\, (\ws_\delta\cap B)={\dim_{\mathcal{H}}}  \, \ws_\de  $ for any non-trivial closed ball $B\subset \zu^d$.  In particular, in general there is no Hausdorff measure whose restriction to $\widetilde {\mathcal{L}} (\Omega,f)$ is positive. 

Let us   illustrate our purpose. 

In $\zuo$, consider the system $\mathcal{R}$ associated with the rational numbers and the function  $f(x)= 1+x$ (the crucial property is the strict monotonicity of $f$). We are interested in $\, \mathcal{L}(\zuo,f)=\{x\in \zuo: \de_x\geq 1+x\}$ and  $\ws(\zuo,f)=\{x\in \zu: \de_x=1+x\}$. 

Jarnik's theorem obviously implies that ${\dim_{\mathcal{H}}} \, \mathcal{L}(\zuo,f) =1$. Indeed, using that $1+x$ tends to $ 1$ when $x>0$ tends to $0$,  for every $\ep>0$, the set $\, \mathcal{L}(\zuo,f)\cap [0,\ep] $ contains all the real numbers whose approximation rate $\de_x$  is larger than $1+\ep$. These real numbers form a set of Hausdorff dimension $1/(1+\ep)$. Letting $\ep$ tend to zero yields the result.
 
Similar arguments imply that for every $\ep>0$,   $\ws(\zuo,f)\cap [\ep,1] $ has  Hausdorff dimension less than $1/(1+\ep)$. However,  Theorem \ref{main} claims that the  Hausdorff dimension of $\ws(\zuo,f)$ is 1. Consequently, the elements of $\ws(\zuo,f)$ responsible for the value of the Hausdorff dimension of $\ws(\zuo,f)$ are  ``localized" around 0. Observe that in this context    $f$ does not reach its infimum.  

For this reason, we refer to Theorem \ref{main}  as a  localized Diophantine approximation. The proof will consist in  constructing a family  of Cantor sets $(\mathcal{K}_\ep)_{\ep>0}$, all included in $ \ws(\zuo^d,f)$, which will be  located closer and closer to one infimum of the function $f$. These Cantor sets will contain elements $x$ with   prescribed approximation rates (which may  depend on $x$). The sequence of dimensions $\dim_{\mathcal{H}} \, \mathcal{K}_\ep$ will be increasing to the desired dimension $\displaystyle \frac{d}{\min \{f(x): x\in \Omega\}}$, as $\ep$ tends to zero.

\ms

%

We  will  prove Theorem \ref{main2}, which is slightly more general than Theorem \ref{main}. This second version is determinant  for its application to the analysis  of the  H\"older singularities of some Markov processes \cite{MARK}. This extension addresses functions $f$ which are continuous outside a set $E$ with a  given  Hausdorff dimension.
\begin{theorem}
\label{main2}
Suppose that  the assumptions of Theorem \ref{main} on the system $\mathcal{S}$ are satisfied:  (\ref{optim}), $\mathcal{C}_1$ and $\mathcal{C}_2$ hold for the system $\mathcal{S}$.

Let  $\Omega$ be a non-empty compact subset of $(0,1)^d$, such that $\overline{\overset{\circ}{\Omega}} =\Omega$. 

Let $E\subset \Omega$ be a subset of $\Omega$.

Let  $f:(0,1)^d\to [1,\infty)$ be a function whose restriction to $\Omega\setminus E$ is continuous. 

 Suppose that $\displaystyle {\dim_{\mathcal{H}}}\, E< \frac{d}{\inf \{f(x): x\in \Omega\setminus E\}}$. Then
\begin{equation}
\label{resultat}
{\dim_{\mathcal{H}}} \, \mathcal{L}(\Omega\setminus  E,f) = {\dim_{\mathcal{H}}} \,  \widetilde { \mathcal{L}}(\Omega\setminus E,f) = \frac{d}{\inf \{f(x): x\in \Omega\setminus E\}}.
\end{equation}
If, moreover, $\displaystyle {\dim_{\mathcal{H}}}\, E< \frac{d}{\sup \{f(x): x\in \Omega\setminus E\}}$, then the sets $\, \mathcal{L}(\Omega\setminus  E,f)$ and $\widetilde { \mathcal{L}}(\Omega\setminus E,f) $ are dense in $\Omega$.

\end{theorem}

In general the sets $\, \mathcal{L}(\Omega\setminus  E,f)$ and $\ws(\Omega\setminus  E,f)$ cannot be studied by Khintchine-like formulas or by mass transference formulas as stated in \cite{Beres,Beres2} (unless a localized version of these theories is developed). Moreover, they do not possess any large intersection properties \cite{FalcLI}, due to the presence of the non-constant function $f$. 


\ms

The paper is organized as follows. Conditions $\mathcal{C}_1$ and $\mathcal{C}_2$, as well as some preliminary results, are given in Section \ref{seccond}. The lower bound in the two-sided equality (\ref{resultat}) is proved in Section \ref{secproof}, while the corresponding upper bound is demonstrated in Section  \ref{secmaj}.  Finally, several examples of suitable systems (including the rational system) are studied in Section  \ref{secexa}.

\section{Definitions and conditions $\mathcal{C}_1$ and $\mathcal{C}_2$}
\label{seccond}
In $\R^d$, we work with the $L^\infty$ norm.

\subsection{Hausdorff measure, gauge functions and Hausdorff dimension}\label{Hausdorff}

Let $\zeta$ be a {\em gauge} function, i.e. a non-negative
non-decreasing function on $\mathbb{R}_+$ such that $\lim _{x\to
0^+}\zeta(x)=0$.  Let $S$ be a subset of $\mathbb{R}^d$. For all
$\eta>0$, let us define the quantity 
\[\mathcal{H}_\eta^\zeta(S) = \inf \ \sum_{i\in \mathcal{I}}
\zeta \left (\left| C_i \right|\right),
\]
the infimum being taken over all the countable families $\{C_i\}_{i\in \mathcal{I}} $ of subsets of $\mathbb{R}^d$ such that $\bigcup_{i\in I} C_i$ is a covering of $S$ and $|C_i| \leq \eta$ for all $i\in I$ ($|C_i|$ stands for the diameter of $C_i$). As $\eta$ decreases to 0,
$\mathcal{H}_\eta^\zeta(S)$ is non-decreasing, and $\mathcal {H}^\zeta(S)=\lim _{\eta\to 0} \mathcal{H}_\eta^\zeta(S)$ defines an outer measure on $\mathbb{R}^d$, called Hausdorff $\zeta$-measure, whose restriction to the Borel  $\sigma$-field is a measure.

Let $\alpha> 0$ be a real number. The $\alpha$-dimensional Hausdorff measure on $\mathbb{R}^d$ is the measure $\mathcal {H}^{\zeta_\alpha}$, where $\zeta_\alpha(x)=x^\alpha$. For each $S\in \mathbb{R}^d$, there exists a unique real number $0\le D\leq d$ such that   $\displaystyle D= \sup\left\{\alpha > 0:\mathcal {H}^{\zeta_\alpha}(S) = +\infty\right\}
= \inf\left\{\alpha\ge 0:\mathcal {H}^{\zeta_\alpha}(S) =0 \right\}$ (with
the convention $\sup \, \emptyset =0$). This real number $D$ is called the Hausdorff
dimension of $S$ and denoted ${\dim_{\mathcal{H}}}\, S$.
We refer the reader to \cite{F2,MATTILA} for more details.

If $m$ is a Borel probability measure over $[0,1]^d$,  then its lower and upper Hausdorff dimensions are respectively defined by  
\begin{eqnarray*}
{\dim_{\mathcal{H}}}_* (m) & = & \inf\{{\dim_{\mathcal{H}}}\, B:m(B)>0\}\\
{\dim_{\mathcal{H}}}^*(m) & = & \inf\{{\dim_{\mathcal{H}}}\, B: m(B)=1\}.
\end{eqnarray*}
When ${\dim_{\mathcal{H}}}_* (m)={\dim_{\mathcal{H}}}^*(m)$, this common value is called the Hausdorff dimension of $m$ and denoted ${\dim_{\mathcal{H}}}(m)$.

\subsection{Notations}

In the rest of the paper, we consider $(x_n)_{n\ge 1}$  a sequence of elements of $[0,1]^d$, and  $(r_n)_{n\ge 1}$ a non-increasing sequence of positive real numbers converging to 0 when $n$ tends to infinity. We then define the {\em system} $\mathcal{S}$ as the sequence of couples $\big((x_n,r_n)\big)_{n\geq 1}$.

\ms

For every  integer $j\ge 0, $ we denote by $\mathcal{G}_j$ the set of dyadic sub-cubes of $[0,1]^d$ of generation $j$, and let $\mathcal{G}_*$ stand for $\bigcup_{j\ge 1}\mathcal{G}_j$.  For any dyadic cube $I\in \mathcal{G}_*$, we set $g(I)=-\log_2(|I|)$, the dyadic generation of $I$ (recall that $|I|$ stands for the diameter of $I$).

\sk

We denote by $\Phi$ the set of   functions $\varphi:\mathbb{R}_+ \to \R_+$  satisfying
\begin{itemize}
\item
$\varphi$ is a non-decreasing continuous functions such that $\varphi(0)=0$,
\item
 $r\mapsto r^{-\varphi(r)}$ is decreasing and tends to infinity as $x>0$ tends to 0, 
 \item
 for all real numbers $\alpha,\beta>0$,    the mapping $r\mapsto r^{\alpha-\beta \varphi(r)}$ is increasing in a   neighborhood of 0. 
\end{itemize}

\ms

We introduce  now   the conditions on the system $\mathcal{S}$. These conditions essentially ensure  an homogeneous repartition in $\zu^d$ of the points   $(x_n)_{n\ge 1}$, and limit the overlaps between the balls $B(x_n,r_n)$.

\subsection{Condition $\mathcal{C}_1$: Weak redundancy}\label{WRe}

\begin{definition} Given the system $\mathcal{S}= \big ((x_n,r_n)\big )_{n\ge 1}$,  we define the irreducible   sub-system   $\big ((y_n,\rho_n)\big )_{n\ge 1}$ associated with $\big ((x_n,r_n)\big )_{n\ge 1}$ as follows:
$$
\big ((y_n,\rho_n)\big )_{n\ge 1}=\big ((x_n,r_n)\big )_{n\ge 1,\ n=\min\{p\ge 1: \, x_p=x_n\}}.
$$
\end{definition}
If $x\in \{x_n:n\ge 1\}$, then  the irreducible subsystem   $\big ((y_n,\rho_n)\big )_{n\ge 1}$ contains  one (and only one) couple of the form $(x,r)$, where $r= \max\{r_n:  (x_n,r_n)\in \mathcal{S}\}$. This definition is needed since the initial system $\big ((x_n,r_n)\big )_{n\ge 1}$  may be very redundant (this occurs when  one element $x$ appears infinitely many times in the sequence $(x_n)_{n\geq 1}$, as in the case of the system of  rational numbers  $\big((p/q, 1/q^2)\big)_{q\geq 1, \, 0\leq p \leq q-1 }$.)

\begin{definition} Let $\big ((x_n,r_n)\big )_{n\ge 1}$ be a system, and consider its irreducible subsystem  $\big ((y_n,\rho_n)\big )_{n\ge 1}$. For any integer $j\ge 0$ we set
\begin{equation}
\label{deftj}
 \mathcal{T}_{j} =\big \{n: \ \ 2^{-(j+1)}<\rho_n\le
 2^{-j}\big \}.
\end{equation}
\end{definition}

 \noindent
 {\bf Condition $\mathcal{C}_1$:} 
The system  $\big ((x_n,r_n)\big )_{n\ge 1}$ satisfies  $\mathcal{C}_1$ when  
 there exists  a non-decreasing sequence of
integers $(N_{j})_{j\ge 0}$ such that
\begin{enumerate}
\item
we have $ \displaystyle\lim_{j\to\infty} \frac{\log_2 N_{j}}{j}=0$.
\item
 for every $j\ge 1$, $\mathcal{T}_j$ can be decomposed into at most
$N_{j}$ pairwise disjoint subsets (denoted
$\mathcal{T}_{j,1},\dots,\mathcal{T}_{j,N_{j}} $) such that for each
$1\le i\le N_{j} $, the balls $B(y_n,\rho_n)$, where $ n$ ranges over $ \mathcal{T}_{j,i}$, are pairwise disjoint.
\end{enumerate}


\ms

Each $\mathcal{T}_{j,i}$ has   cardinality less  than $2^{d(j+1)}$, and $\mathcal{T}_j$ has   cardinality less than $N_j \cdot 2^{d(j+1)}$. 

{Condition $\mathcal{C}_1$}  ensures that every $t\in\zu^d$ is covered by at most $N_j$ balls of the form $B(y_n,\rho_n)$, $n\in \mathcal{T}_j$. The fact that  $N_j$ does not increase too fast toward infinity explains the appellation  ``weak redundancy'' given to $\mathcal{C}_1$ in \cite{BSUbiquity}. 

\subsection{Condition  $\mathcal{C}_2$: a fine non-overlapping condition}

In order to obtain Theorems \ref{main} and \ref{main2}, an additional  property is required on the  system.
We emphasize that   $\mathcal{C}_2$, though technical,  is satisfied by many natural systems, as explained in Section \ref{secexa}. It appears that, except for the Poisson system,  $\mathcal{C}_1$ and $\mathcal{C}_2$ are quite easy to check.

\begin{definition}
\label{def3}
Suppose that  $\mcs=\big ((x_n,r_n)\big )_{n\ge 1}$ satisfies  $\mathcal{C}_1$, and consider the sequence $(N_j)_{j\geq 1}$ associated  with $\mcs$ by $\mathcal{C}_1$.

There exists a continuous function $\psi:\R_+\to\R_+$ such that $\psi(0)=0$ and for every $j\geq 1$, $N_j$ can be written as 
\begin{equation}
\label{defpsi}
N_j=2^{dj\psi(2^{-j})}.
\end{equation}
For every  $\varphi\in\Phi$  and for every $j\ge 1$, we define
\begin{eqnarray}
\label{defgamma}
\gamma (j) & = & \max \left\{k\in\mathbb{N}: N_k2^{dk}\le 2^{-dj\varphi(2^{-j})}2^{dj}\right\}\\
\nonumber& =&\max \left\{k\in\mathbb{N}:  2^{dk(1+\psi(2^{-k}))}\le 2^{dj(1-\varphi(2^{-j}))} \right\}.
\end{eqnarray}
\end{definition}
Obviously $\gamma(j)\leq j$, and the difference $j-\gamma(j)$ can be written as  $j \theta(2^{-j})$, where the mapping $\theta:\R_+\to\R_+$ is continuous and $\theta(0)=0$. 

The sequences  $(\gamma (j))_{j\ge 1}$ and $(\theta(2^{-j}))_{j\ge  1}$ depend on the sequence $(N_j)_{j\ge 1}$ and on $\ph$.  Nevertheless, in the following, we omit to write this dependence, since by Property $\mathcal{C}_2$, both $(N_j)_{j\ge 1}$ and  $\ph$ will be fixed once for all.

\begin{definition}\label{def2}
Let $\varphi \in \Phi$ and  $(N_j)_{j\ge 1}$ be defined as in Definition \ref{def3}. Let  $V\in \mathcal{G}_*$ be a dyadic cube in $\zu^d$. Let $\delta>1$ be a real number.

Recall that $g(V)= -\log_2 |V|$ if the dyadic generation of $V$.

The property $\mathcal{P}(V,\delta )$ is said to hold   when  there exists $x(V)\in V\subset \zu^d$ and a positive real number $r(V)$  verifying:
\begin{itemize}
\item
 $(x(V),r(V))\in \mcs$,
\item
 $2^{-g(V) -1} \leq  r(V) < 2^{-g(V) }$,
\item and 
\begin{eqnarray*}
\!\!\!\!\! \!\!\!\!\!\!\! &&B\big(x(V),r(V)^{\delta} \big) \,  \bigcap \, \Big\{   x_p :    \gamma\big (g(V)\big ) \le   - \log _2 r_p      <   [\delta (g(V)+1)]+4 \Big\} \\ 
\!\!\!\!\! \!\!\!\!\!\!\!&& =\{x(V)\}.
\end{eqnarray*}
\end{itemize}
\end{definition}
The notation $[y]$ stands for the integer part of the real number $y$.

Recall that $\gamma\big (g(V) \big )  \leq   g(V)$,  and note that $ [\delta g(V)]$ is heuristically  the generation of the largest dyadic cube  included in   the contracted ball $
B\left(x(V),r(V)^{\delta} \right) $.

\begin{center}
\begin{figure}
\begin{center}\label{fig.1}
Figures are available on our web site.
\caption{Property $\mathcal{P}(V,\delta)$ }
\end{center}\end{figure}
\end{center}

$\mathcal{P}(V,\delta)$  holds   when, except $x(V)$,   all the elements $x_p$, where $p$ ranges over the indices such that  $ \gamma\big (g(V)\big )   \le - \log _2 r_p  <  [\delta (g(V)+1)]$ $+4$,  avoids the contracted ball $
B\left(x(V), r(V) ^{\delta} \right) $  (see Figure 2.).  The constant 4 is due to  technicalities along  the proof. Note that $\mathcal{P}(V,\de)$ depends on $(N_j)_{j\geq 1}$ and $\ph$  via  $\gamma$ (formula (\ref{defgamma})), but as said above we do not mention this dependence since $(N_j)_{j\geq 1}$ and $\varphi$ are fixed by $\mathcal{C}_2$.

\ms 

$\mathcal{P}(V,\delta )$ seems to be a reasonable property, maybe not for all   dyadic cubes $V$, but at least for a large number among them.  {Condition  $\mathcal{C}_2$} is meant to ensure the validity of $\mathcal{P}(V,\delta )$ for a sufficient set of cubes $V$ and approximation rates $\delta$.

\medskip

{\bf Condition  $\mathcal{C}_2$:} A system $\mcs$ satisfies  $\mathcal{C}_2$ when $\mcs$ satisfies  $\mathcal{C}_1$ and when there exists :
\begin{itemize}
\item
a function $\varphi\in\Phi$, 
\item
a non-decreasing sequence of integers   $(N_j)_{j\ge 1}$ as in Definition \ref{def3},
\item
a continuous function $\kappa:(1,+\infty) \to (0,1]$,
\item
a dense subset $\Delta$ of $(1,\infty)$,
\end{itemize}
with the following property: \\
For every $\delta\in \Delta$,  for every dyadic  cube $U$ of $\zu^d$, there are infinitely many integers $j\ge g(U)$ satisfying
\begin{equation}
\label{minqjd}
\#\mathcal{Q}(U,j,\delta) \, \ge \, \kappa(\delta )  \cdot  2^{d(j -g(U))} ,
\end{equation}
where
$$
\mathcal{Q}(U,j,\delta)=\left\{V\in \mathcal{G}_j:\ V\subset U \ \mbox{ and } \ \ \mathcal{P}\big (V,\delta \big )\text{ holds}\right \}.
$$

\ms

In the following, the system $\mcs$ satisfies ($\mathcal{C}_1$ and) $\mathcal{C}_2$. Hence, $\ph$ and $(N_j)_{j\ge 1}$  are given, and all the parameters introduced  from now on depend on them.

Observe  that $2^{d(j-g(U))}$ is the number of dyadic cubes $V$ of generation $j\geq g(U)$ included in $U$. Among them, $\mathcal{Q}(U,j,\delta)$ contains the cubes enjoying the property   $\mathcal{P}\big (V,\delta \big )$. As claimed above, {Condition  $\mathcal{C}_2$} guarantees that  given a dyadic cube $U$ and $\de\in\Delta$, infinitely often   a given proportion of the dyadic subcubes $V$ of generation $j$ included in   $U$ satisfies $ \mathcal{P}\big (V,\delta \big )$.



\begin{remark}
For the rational system and other deterministic systems provided in Section \ref{secexa}, the function $\kappa$ can be taken constant: $\forall \, \de>1$, $\kappa(\de)=\kappa \in (0,1)$. The possible dependence in $\delta$ of the factor $\kappa$ is introduced to include the systems obtained as Poisson point processes in the upper-half plane. This is explained in Section \ref{secexa}.
\end{remark}

\subsection{A preliminary result }

We shall need the following lemma, which requires only $\mathcal{C}_1$.

\begin{lemma}\label{control} Let $\big ((x_n,r_n)\big )_{n\ge 1}$ be a system  and  let    $\big ((y_n,\rho_n)\big )_{n\ge 1}$ be the corresponding irreducible subsystem. 
Suppose that $\mathcal{C}_1$ is satisfied. 

For every $\de>1$, for every dyadic cube  $U\in\mathcal{G}_*$, and  every integer $j\ge \delta \cdot g(U)$, let us introduce the set of cubes
$$
\widetilde {\mathcal{Q}}(U,j,\delta )\!=\!\left\{ \! V\in \mathcal{G}_j:V\subset U , \,   V \ \bigcap \left( 
\bigcup_{k=g(U)}^{\gamma(j)} \ \bigcup_{p\in \mathcal{T}_k}B(y_p, (\rho_p)^\delta) \!\right )\neq\emptyset\right\}.
$$
Then, there exists a constant $C_d$ depending only on $d$ such that
\begin{equation}
\label{eq1} 
\# \widetilde{\mathcal{Q}}(U,j,\delta )\!  \leq C_d\cdot 2^{d(j-g(U))} \!\left[   2^{-dj\varphi(2^{-j})}+ \!\! \sum_{g(U)\le k\le j/\delta}  2^{-dk(\delta-1-\psi(2^{-k}))}
\right].
\end{equation}
\end{lemma}

The sets $ \widetilde{\mathcal{Q}}(U,j,\delta )$ contains the dyadic cubes of generation $j$ which  intersect the irreducible balls $B(y_n,\rho_n)$ when $n$ ranges in $\mathcal{T}_p$, $p\in [g(U),g(U)+1,...,\gamma(j)]$. It is crucial  in the further construction of Cantor sets that $ \widetilde{\mathcal{Q}}(U,j,\delta )$ cannot contain a too large number of cubes  (see Figure 3.).

\begin{center}
\begin{figure}
\begin{center}\label{fig2}
Figures are available on our web site.
\caption{Property $\widetilde{\mathcal{Q}}(U,j,\delta)$ }
\end{center}\end{figure}
\end{center}

Our proof will use the following standard estimates.
\begin{lemma}\label{intersection} There exists a constant $C'_d$ depending on $d$ only such that:
\begin{enumerate}
\item If $r_0>0$, $j\in\N$ and $B$ is a closed ball such that $2^{-j}\le |B|\le r_0$, then $B$ intersects at most $C'_d\cdot 2^{dj} \cdot  r_0$ elements of $\mathcal{G}_j$.
\smallskip
\item If $U\in\mathcal{G}_*$, $k$ is    an integer larger than $g(U)$ and $\mathcal{T}$ is a family of pairwise disjoint closed balls of radius larger than $2^{-(k+1)}$, then $U$ intersects at most $C'_d \cdot 2^{d(k-g(U))}$ elements of  $\mathcal{T}$.
\end{enumerate}
\end{lemma}

\noindent{\it Proof of Lemma \ref{control}.} 
Let $U\in\mathcal{G}_*$ and $j\ge \delta \cdot g(U)$. Let $k$ be an integer such that   $k\in [g(U),\gamma (j)]$. 

We are going to count the number of dyadic cubes  $V$ in $\mathcal{G}_j$  which are  included in $U$ and which intersect  balls of the form $B(y_p, (\rho_p)^\delta)$ for some  $p\in\mathcal{T}_k$.  Two cases shall be distinguished:

\sk

$\bullet$  If $g(U)\leq k\le j/\delta$, then for $p\in\mathcal{T}_k$ we have $2^{-j}\le 2^{-k\de } \le |B(y_p, (\rho_p)^\delta )|\le 2^{1-k\delta}$. Consequently, $B(y_p, (\rho_p)^\delta )$ intersects at most $C'_d\cdot 2^{dj}2^{d(1-k\delta)}$ elements of $\mathcal{G}_j$.

Moreover,  by construction,  $\mathcal{T}_k = \bigcup_{1\le l\le N_k} \mathcal{T}_{k,l}$, where the elements of each $\mathcal{T}_{k,l}$ are pairwise disjoint closed balls of radius larger than $2^{-(k+1)}$. Consequently, for  $1\le l\le N_k$, the cardinality of those  integers $p \in \mathcal{T}_{k,l}$ satisfying  $B(y_p, \rho_p )\cap U \neq \emptyset$ is at most $C'_d\cdot 2^{d(k-g(U))}$. Thus, there are at most $C'_d\cdot N_k2^{d(k-g(U))} $ integers $p\in \mathcal{T}_k$ satisfying  $B(y_p, \rho_p )\cap U \neq \emptyset$.

Combining the last  remarks,  the  cardinality of the subset of $\mathcal{G}_j$  whose elements are included in $U$ and meet a ball $B(y_p, (\rho_p)^\delta)$ with $p\in\mathcal{T}_k$, is less than 
\begin{eqnarray*} 
(C'_d)^2 N_k 2^{d(k-g(U))} \cdot 2^{dj}2^{d(1-k\delta)}& = & 2^d(C'_d)^2  \cdot 2^{d(j-g(U))} N_{k} \cdot  2^{-dk(\delta-1)}.
 \end{eqnarray*} 

\sk

$\bullet$  
If $j/\delta<k\le \gamma (j)$, then  for every $p\in\mathcal{T}_k$ we have $|B(y_p, (\rho_p)^\de)|\le 2^{1-j}$.
Hence $B(y_p, (\rho_p)^\de)$ intersects at most $3^d$ cubes of $\mathcal{G}_j$. Consequently, the  cardinality of the subset of $\mathcal{G}_j$  whose elements are included in $U$ and meet a ball $B(y_p, (\rho_p)^\delta)$ with $p\in\mathcal{T}_k$ is at most $3^dC'_d \cdot  N_k2^{d(k-g(U))}$.

\ms

Summarizing the above estimates, we obtain
\begin{eqnarray*} 
\# \widetilde{\mathcal{Q}}(U,j,\delta )&\le &2^d(C'_d)^2   \cdot 2^{d(j-g(U))}  \sum_{g(U)\le k\le j/\delta} N_k2^{-kd(\delta-1)}
\\ && + \,3^d C'_d \cdot 2^{-dg(U)}\sum_{j/\delta<k\le\gamma (j)} \, {N}_{k}2^{dk}.
\end{eqnarray*}
Using that   $( N_k)_{k\ge 1}$   is non-decreasing, we get
\begin{eqnarray*} 
3^dC'_d\cdot 2^{-dg(U)}\sum_{j/\delta<k\le\gamma (j)]} \,{N}_{k}2^{dk}
&\le &   3^dC'_d\cdot 2^{-dg(U)} {N}_{\gamma (j)}\sum_{j/\delta<k\le\gamma (j)} \, 2^{dk} \\
&\leq &2\cdot 3^dC'_d\cdot N_{\gamma (j)}\, 2^{d\gamma (j)} \, 2^{-dg(U)} \\ & \leq & 2\cdot 3^dC'_d\cdot 2^{d(j-g(U))}  2^{-dj \varphi(2^{-j})} ,
\end{eqnarray*}
where we used the definition of $\gamma (j)$ in the last inequality.
Moreover, using the definition (\ref{defpsi}) of $\psi(2^{-k})$ based on  $N_k$, we find that
$$
2^d(C'_d)^2  \cdot 2^{d(j-g(U))}  N_{k}  2^{-kd(\delta-1)}
\leq  2^d(C'_d)^2 \cdot 2^{d(j-g(U))}    2^{-kd(\delta-1-\psi(2^{-k}))}.$$
Equation (\ref{eq1})   follows easily. 

\section{Lower bound for the Hausdorff dimensions in Theorem \ref{main2}}
\label{secproof} 

Let  $\mcs=\big((x_n,r_n)
\big)_{n\geq 1}$ be a system satisfying $\mathcal{C}_1$, and let $\varphi$, $(N_j)_{j\geq 1}$, and $\Delta$ be fixed so that $\mcs$ satisfies also   $\mathcal{C}_2$. We denote by  $\big((y_n,\rho_n)
\big)_{n\geq 1}$ the irreducible subsystem of $\mcs$.

\ms
Consider $\Omega$, $E$ and $f$ as in  Theorem \ref{main2}. We have 
 \begin{equation}
 \label{defh}
 h:=\frac{d}{\inf\{f(x): x\in \Omega \setminus E\}}>{\dim_{\mathcal{H}}}\, E.
 \end{equation}
Our aim is to prove that the Hausdorff   dimension  of $\ws(\Omega \setminus E,f)=\{x\in \Omega\setminus E:\delta_x=f(x)\}$ equals $h$. 
  We are going to construct a family of Cantor sets all included in $\ws(\Omega \setminus E,f)$  and such that the supremum of their Hausdorff dimensions is larger than (or equal to) $h$.

\subsection{First simplifications}

Before starting the constructions, we make some remarks:

\begin{itemize}
\item
If the restriction of $f$ to $\Omega\setminus E$ is equal to the minimum of $f$ over $U \cap{\overset{\circ}{\Omega}} \setminus E$, where $U$ is a non-empty dyadic cube, then the result follows from Theorem \ref{thubi}. Thus we will assume that this is not the case, i.e. the subset of $\Omega\setminus E$ over which $f$ reaches its minimum is nowhere dense (this set is empty in general).
\item
Problems may occur in the construction below when $h=d$, i.e. when $\inf \{f(x): x\in \Omega\setminus E\} =1$. Let us explain how we circumvent such difficulties. Assume that $\inf \{f(x): x\in \Omega\setminus E\} =1$. For every $\ep>0$ small enough, it is possible to find  a dyadic cube $U_\ep$ such that the restriction of $f$ to $U_\ep \cap( \Omega\setminus E)$ has an infimum  which belongs to the open interval $(1,1+\ep)$. The construction below can be applied to the set $\ws(U_\ep \cap(\Omega \setminus E),f)$, and we find that $  \dim_{\mathcal{H}} \ws(U_\ep \cap(\Omega \setminus E),f) \geq d/(1+\ep)$. Letting $\ep$ tend to zero yields the result. 
\end{itemize}

Thus we assume that   $h$ defined by  (\ref{defh}) is strictly less than $d$.

\subsection{Preliminary work}

Fix $\varepsilon\in (0,  h)$, and recall that $h<d$. By definition (\ref{defh}) of $h$, and since ${\dim_{\mathcal{H}}} \, E<d={\dim_{\mathcal{H}}} \,  \Omega$, there exists $y_\ep\in {\overset{\circ}{\Omega}} \setminus E$ such that $\displaystyle h-\ep/2 \leq \frac{d}{f(y_\ep)} \leq h $. Hence, using the continuity of $f$ at $y_\ep$,  for every $y $ in a neighborhood $\Omega_\ep \subset  {\overset{\circ}{\Omega}} $ small enough around $y_\ep$,  we have $ \displaystyle h -\ep\le \frac{d}{f(y)}\le h$. Equivalently, when $\ep$  small enough, we have
\begin{eqnarray}
\label{encad}
\forall \ y\in\Omega_\ep, \ \ \frac{d}{h} \leq f(y) \leq \frac{d}{h-\ep} \leq \frac{d}{h} \left( 1 + 2 \frac{\ep}{h} \right).
\end{eqnarray}

\ms

Recall that $\Delta$ is the set of admissible approximation rates allowed by property $\mathcal{C}_2$. In every dyadic cube $V\in \mathcal{G}_*$    included in $ \Omega _\ep$, we pick up  an element $y_V\in V$ and we choose a real number  $\delta(V)\in \Delta$ such that ($\ph$ is fixed by $\mathcal{C}_2$ and $\psi$ is defined by (\ref{defpsi}))
\begin{equation}
\label{defdeltau}
\delta(V)\in \big [\displaystyle f(y_V) + d(\ph(|V|) +\psi(|V|))  , f(y_V) +3d(\ph(|V|) +\psi(|V|))  \big ].
\end{equation}

Observe that the real numbers $\delta(V)$ are bounded from above and below, since $\phi$ and $\psi$ are continuous and $f$ is bounded on $\Omega_\ep$.  Moreover,   by formula (\ref{encad}) there exists a constant $\alpha>1$ such that for every $V$ having diameter small enough  one has
\begin{equation}
\label{minde}
\delta(V)  -3d\ph(|V|) \geq \alpha >1.
\end{equation} 
 Since the function $\kappa(\cdot)$ determined by condition $\mathcal{C}_2$ is continuous,  there is a constant $\kappa \in (0,1)$ such that  for every $\de $ belonging to the set $ \{\delta(V):     V \mbox{ dyadic cube } \subset \Omega_\ep\}$,    for every   dyadic cube $U \subset \Omega_\ep$, (\ref{minqjd}) holds infinitely often  with the same constant $\kappa $ (instead of $\kappa(\de)$). We   choose $\kappa$   so that $2^{d+1}/\kappa$ is a positive power of $2$. This will simplify a little bit the forthcoming  constructions.

\ms\ms

We now start the construction of a Cantor set $\mathscr{K}_\varepsilon$ such that $\mathscr{K}_\varepsilon\setminus E$ is included in $\ws(\Omega \setminus E,f) \cap \Omega_\ep$ and simultaneously a probability measure $\mu_\varepsilon$ supported by  $\mathscr{K}_\varepsilon$ such that $\displaystyle{\dim_{\mathcal{H}}}_*( \mu_\varepsilon)\ge h_\varepsilon$, for some real number $h_\ep$ which satisfies  $\lim_{\ep\to 0^+}h_\ep=h$. 

Assume for a while that the construction of  $\mu_\varepsilon$ and   $\mathscr{K}_\varepsilon$ is achieved. 
Then the  lower bound in (\ref{resultat}) of Theorem \ref{main2} is obtained by the following argument. Since ${\dim_{\mathcal{H}}}\,  E<h$, when $\ep$ is small enough we have $\mu_\varepsilon(E)=0$. Recalling that the support of $\mu_\varepsilon$ is $\mathscr{K}_\varepsilon$,  we deduce that
 $$\displaystyle {\dim_{\mathcal{H}}} \, \ws(\Omega \setminus  E,f)\geq  {\dim_{\mathcal{H}}}\,\mathscr{K}_\varepsilon\setminus E \geq {\dim_{\mathcal{H}}} ( \mu_\varepsilon)   \ge h_\ep.$$
 Letting $\ep$ tend to $0$ yields $$\displaystyle {\dim_{\mathcal{H}}} \, \ws(\Omega \setminus E,f)\geq  h.$$

\ms

The Cantor set $\mathscr{K}_\varepsilon$ will be obtained as a limsup set of the form
$$\mathscr{K}_\varepsilon = \bigcap_{n\ge 0} \ \ \bigcup_{U\in\mathcal{F}_n}  \ \   U ,$$
where for every $n\ge 0$, $\mathcal{F}_n$ is a collection of pairwise disjoint closed dyadic cubes $U$ such that each element of $\mathcal{F}_{n+1}$ is included in one (and by construction only one)  element of $\mathcal{F}_n$.

\ms

The sequence $(\mathcal{F}_n)_{n\ge 0}$ is built by induction, as follows.

\ms
 
 At first, we choose a dyadic cube $U_0$ included in $\Omega_\ep$, small enough so that    $3 d(\ph(  |U_{0}| ) + \psi(|U_{0}|)  ) \le   \varepsilon $ and $ \displaystyle \frac{2^{d+1}}{  \kappa}\le |U_0|^{-d\varphi(|U_0|)}$, where we recall that $\kappa$ is the constant appearing in (\ref{minqjd}) (the dependence on $\de$ has been removed by an argument above). We define   $\mathcal{F}_{0}=\{U_{0}\}$. This choice also implies that for any dyadic cube $U\subset U_0 \subset \Omega_\ep$ that we are going to consider, we have (using (\ref{encad}) and (\ref{defdeltau}))
\begin{equation}
\label{eq00}
 \delta(U) \leq f(x_\ep)+ \ep \left(2 \frac{d}{h^2} +1\right)=: H_\ep.
 \end{equation}

%


\subsection{Construction of the first generation of the Cantor set, $\mathcal{F}_1$} 

Let us find the elements of $\mathcal{F}_1$.

We apply property $\mathcal{C}_2$ and Lemma~\ref{control} with $U=U_0$ and $\delta=\delta(U_0)$.  This yields that there are infinitely many integers  $j\ge g(U_0)$ such that 
\begin{eqnarray*}
\#\mathcal{Q}(U_0,j,\delta(U_{0})) & \ge &  \kappa  \cdot  2^{d(j-g(U_0))}\\
\mbox{and }\#\widetilde {\mathcal{Q}}(U_0,j,\delta(U_{0}))  & \le &    C_d  \cdot 2^{d(j-g(U_0))}  \Big( 2^{-dj\varphi(2^{-j})}+\\ && \sum_{g(U_0)\le k\le j/\delta(U_{0})}    2^{-dk(\delta(U_{0})-1-\psi(2^{-k}))}
\Big).
\end{eqnarray*}
We use (\ref{minde}) to  bound from above the sum in the second equation:
\begin{eqnarray*}\sum_{g(U_0)\le k\le j/\delta(U_{0})}   \!  \!  \!  \!  \! 2^{-dk(\delta(U_{0})-1-\psi(2^{-k}))}  & \leq &   \sum_{g(U_0)\le k }    2^{-dk \alpha}   \leq C  \cdot 2^{-d \alpha g(U_0) }.
\end{eqnarray*}
The constant $C$ does not depend on $U_0$. Consequently, the second upper bound above can be simplified into
\begin{eqnarray*}
\#\widetilde {\mathcal{Q}}(U_0,j,\delta(U_{0})) & \le & C_d\cdot C \cdot 2^{d(j-g(U_0))} (2^{ -j\varphi(2^{-j})} +   2^{-d\alpha g(U_0)}).
\end{eqnarray*}
Provided that $U_0$  has diameter small enough and $j$ is large enough, we have $\#\widetilde {\mathcal{Q}}(U_0,j,\delta(U_{0}))  \leq \kappa/4 \cdot 2^{d(j-g(U_0))}$.
From the  inequalities between the cardinalities of   $  {\mathcal{Q}}(U_0,j,\delta(U_{0})) $  and $\widetilde {\mathcal{Q}}(U_0,j,\delta(U_{0})) $, we deduce that for   $j$ large enough, there is  a subset $\widetilde {\mathcal{F}}_1$ of cardinality at least $\displaystyle \frac{\kappa}{2} \cdot  2^{d(j-g(U_0))}$ in $\mathcal{Q}(U_0,j,\delta(U_{0}))\setminus \widetilde {\mathcal{Q}}(U_0,j,\delta(U_{0}))$. Moreover, 
we can find at least $\#\widetilde{\mathcal{F}}_1/2^d$ elements  of $\widetilde{\mathcal{F}}_1$ which are distant from each other by at least $2^{-j}$.   Consequently, we can assume that there are exactly  $\displaystyle \frac{\kappa}{2^{d+1}}  \cdot  2^{d(j-g(U_0))}$  dyadic  cubes  in $\widetilde{\mathcal{F}}_1$, whose mutual distance is  at least $2^{-j}$. 

By construction,  each cube $\widetilde V \in\widetilde{\mathcal{F}}_1$ satisfies simultaneously $\widetilde V\in \mathcal{G}_j$, $\widetilde V\subset U_0$,    $ \mathcal{P}\big (\widetilde V, \delta(U_{0}) \big )   $ and (recall that $j=g(\widetilde V)$)
\begin{eqnarray*}
  \widetilde V \ \bigcap \Big( \
\bigcup_{k=g(U_0)}^{\gamma(j)} \ \ \bigcup_{p\in \mathcal{T}_k} \ B(y_p, (\rho_p)^{\delta(U_{0})}) \Big) = \emptyset.
\end{eqnarray*}
 
 Combining the information, each cube $\widetilde V \in\widetilde{\mathcal{F}}_1$  contains  an element  $x(\widetilde V)$ such that   $( x(\widetilde V),r(\widetilde V)) \in \mcs$  for some radius $r(\widetilde V)$ satisfying        $  2^{-j-1} \le   r(\widetilde V)  < 2 ^{-j} $. By construction we have  
\begin{eqnarray*}
&& B\left(x(\widetilde V),r(\widetilde V)^{{\delta(U_{0})} } \right) \, \bigcap \, \Big\{   x_p : \  \gamma(j)\le- \log _2 r_p  < (j+1) {\delta(U_{0})} +4 \Big \} \\ && =\{x(\widetilde V)\}.
\end{eqnarray*}

\ms

In order to compute the Hausdorff dimension of the limsup sets we are interested in, we must find points with a  prescribed approximation rate in a very  precise way. For this,  a new definition is needed.

\begin{definition}
\label{def1}
For $\delta>1$, $\varepsilon>0$, $x\in \R^d$ and $r>0$ we define the annulus
$$
A(x ,r ,\delta,\varepsilon)=B(x ,r^{\delta })\setminus B(x ,r^{\delta+\varepsilon}).
$$
\end{definition}

For each $\widetilde V\in\widetilde{\mathcal{F}}_1$, consider the associated    annulus 
$$
A(\widetilde V) = A(x(\widetilde V),r(\widetilde V),\delta(U_{0}), \varphi (2^{-j})).$$ 

\begin{remark}
\label{rem1}
The diameter of  $
A(\widetilde V)$  is $2\cdot   r(\widetilde V)^{{\delta(U_{0})} } = 2^{1+\log_2  (r(\widetilde V)^{\delta(U_{0}) })}$. Provided that $j$ is taken large enough, the ``hole" in the annulus $A(\widetilde V)$ is extremely small, since the ratio $ \frac{ r(\widetilde V)^{\delta(U_{0}) }}   {  r(\widetilde V)^{\delta(U_{0})+\varphi (2^{-j})}} =r(\widetilde V)^{ - \varphi (2^{-j})}$  tends to infinity when $j$ tends to infinity (recall that $\ph$ belong to the functional space $\Phi$ and $r(\widetilde V)\sim2^{-j}$).
\end{remark}

Let $ \widetilde{\widetilde V}$ be  one  of the largest closed dyadic cubes included in  $A(\widetilde V)  \bigcap \widetilde V$. 
Using  Remark \ref{rem1}, the generation $g(  \widetilde{\widetilde V})$ of the dyadic  cube $ \widetilde{\widetilde V}$ is at most equal to $[- \log_2  (r(\widetilde V)^{{\delta(U_{0})}} )]+3$.  

{\bf We choose then $V$ to be one of the subcubes of  $ \widetilde{\widetilde V}$ of generation $g(  \widetilde{\widetilde V})+1$ among those cubes of this generation which  are  the closest to $x(\widetilde V)$.} We obtain that:
\begin{itemize}
\item the dyadic generation of $V$ satisfies
 \begin{equation} 
 \label  {encadgv}
\hspace{-4mm} [- \log_2  (r(\widetilde V)^{\delta(U_{0})  } )]\le g(V)\le [- \log_2  (r(\widetilde V)^{\delta(U_{0}) } )]+4,
 \end{equation}
 \item
  for each $p$ such that $\gamma(j)\le -\log_2 r_p \leq [(j+1)\delta(U_{0}) ]+4$, if $x_p\neq x_n=x(\widetilde V)$ then $x_p\not\in B(x(\widetilde V), r(\widetilde V)^{\delta(U_{0}) })$,
  \item   for each $p$ such that $\gamma(j)\le -\log_2 r_p \leq [(j+1)\delta(U_{0}) ]+4$, for all $x\in V$ we have (using the function $\theta$ defined in Definition~\ref{def3})
\begin{eqnarray*}
|x-x_p|\ge |V|&\ge &r(\widetilde V)^{\delta(U_{0}) }/16 \  =  \ 2^{(\gamma(j)-j) \delta(U_{0} )}\cdot 2^{-\gamma (j) \delta(U_{0} )}/16\\
&\ge& 2^{-j\theta(2^{-j}) \delta(U_{0}) } \cdot r_p^{\delta(U_{0})   }/16 \  \ge  \ r_p^{\delta(U_{0}) +\theta (2^{-j})\de(U_0)/2}/16.
\end{eqnarray*}
The last inequality follows from the fact that, when $j$ is large enough, $r_p \leq 2^{-\gamma(j)}\leq 2^{-j/2}$. Using (\ref{eq00}), we see that
\begin{eqnarray*}
|x-x_p|\ \ge  \ r_p^{\delta(U_{0}) +\theta (2^{-j}) H_\ep/2}/16.
\end{eqnarray*}
\end{itemize}
When two dyadic cubes $V$ and $\wv$ are related via such a relationship, we say that {\bf $V$ is the contracted descendant of $\wv$}.

\ms
  
The previous construction guarantees  that (recall that $j=g(\widetilde V)$):
\begin{itemize}
\item (\ref{encadgv}) holds, 
\item
 since   $V\subset  A(\widetilde V)$,  every element $x\in V$ is approximated at a rate  $\in [\delta(U_{0})   , \delta(U_{0})  +\ph(2^{-j})]$ by $(x(\widetilde V),r(\widetilde V))$,
\item
for every $x\in V$, for every $k\in \{g(U_0),...,\gamma (j)\}$, for every $p\in  \mathcal{T}_k  $, $x \notin  B(x_p, (r_p)^{\delta(U_{0}) }) $, i.e. $x$ is not approximated at rate larger than $\delta(U_{0}) $ by these  couples $(x_p,r_p) \in \mcs$ ,

\item
for every $x\in V$, for every integer $p$ such that
$$ \gamma({j}) \, \le \,[- \log _2 r_p ]\, <\,    [- \log_2  (r(\widetilde V)^{\delta(U_{0})  } )] +4,$$  
we have that $x \notin  B\big (x_p, (r_p)^{\delta(U_{0})  +\theta (2^{-j})H_\ep/2}/16\big )$.

\item
The first, third and fourth previous items imply that  if $p$ is such that $2^{-g(V)}\le r_p\le 2^{-g(U_0)}$ and $x_p\neq x(\widetilde V)$, then for all $x\in  V$ we have $x \notin  B\big (x_p, r_p^{\delta(U_{0}) +\theta (2^{-j})H_\ep/2}/16\big )$.
 \end{itemize}


\sk

Since this situation occurs for an infinite number of generations $j$, we choose $j$ large enough so that
\begin{eqnarray}
\label{eq3}
j\ge 2g(U_0) \ \mbox{ and }  \ \max \left( \frac{2^{d+1}} {\kappa} ,  2^{dg(U_0)} \right) \le 2^{dj\varphi(2^{-j})}.
 \end{eqnarray}
The previous inequality  ensures that $\frac{2^{d+1}} {\kappa}\le |V |^{-d\varphi (|V|)}$. This will play a role in Section~\ref{scalingprop}.

By construction, we have $|V|\ge 2^{-(j+1) \delta(U_{0}) }/16$. Consequently, using   (\ref{eq00}) (i.e. $\de(U)$ is   bounded above by $H_\ep$  independently of $U$), without loss of generality one can suppose that  $j$ is large enough so that for some constant $C>0$ (depending on $H_\ep$), 
\begin{eqnarray}\label{eq4}
\mbox{for every }
 \widetilde V\in \mathcal{Q}(U_0,j,\delta(U_{0}) ), \ \ \displaystyle |V|\ge C \cdot {2^{- j \delta(U_{0}) }}.
 \end{eqnarray}
This yields a precise relationship between the diameter of a  cube $\widetilde V\in \widetilde{\mathcal{F}}_1$ and the diameter of its contracted descendant $V$.
 
 \ms

Now, let us  consider the set of contracted descendants of the elements of $\widetilde{\mathcal{F}}_1$
$$
\mathcal{F}_1=\{V: \mbox{$V$ is the contracted descendant of some  }\widetilde V\in \widetilde {\mathcal{F}}_1\}. 
$$

\medskip

We construct a measure $\mu_\varepsilon$ on the algebra $\sigma_1=\sigma (V: \, V\in \mathcal{F}_1)$  generated by the dyadic cubes of $\mathcal{F}_1$  by imposing:
$$
\forall \, V
\in \mathcal{F}_1, \ \ \mu_\varepsilon (V)=(\# \mathcal{F}_1)^{-1}.
$$
Let $ V
\in \mathcal{F}_1$. Recalling that $\# \widetilde{\mathcal{F}}_1=\frac{\kappa}{2^{d+1}} \cdot 2^{d(j-g(U_0))}$, using (\ref{eq3}) we get
$$
\mu_\varepsilon(V)\le \frac{2^{d+1}}{\kappa}  \cdot 2^{d(g(U_0)-j)} \le
 2^{2dj\varphi(2^{-j}) -dj} .$$
Using (\ref{eq4}) we find that  for some universal constant $C>0$ 
$$ 
\mu_\varepsilon(V)\le C \cdot 2^{2dj\varphi(2^{-j})} |V|^{d/ \delta(U_{0})  }.
$$
Due to the monotonicity of $r^{-\varphi(r)}$, we have $2^{2dj\varphi(2^{-j})} \le |V|^{-2d\varphi (|V|)}$, and when $j$ is chosen large enough, $C \leq |V|^{-2d\varphi (|V|)}$. All these computations yield
\begin{equation}
\label{maj1}
\forall\ V\in\mathcal{F}_1,\  \ \ \mu_\varepsilon(V)\le   |V|^{d/\delta(U_{0})  -3d\varphi (|V|)}
\end{equation}


We now fix the integer $j=j_0$ so that all the assumptions above are satisfied. The last property of these cubes of first generation is that for every $V\neq V' \in \mathcal{F}_1$, the distance between $V$ and $V'$ is greater than $2^{-j_0}$.

 
\subsection{Construction of   $\mathcal{F}_2$, the second generation of the Cantor set}

The second generation is obtained as follows.

\medskip

Note that, thanks to (\ref{eq3}),  we insured in the previous step that for each $  U_1\in\mathcal{F}_1$ we have $\frac{2^{d+1}}{\kappa}\le |U_1|^{-d\varphi(|U_1|)}$.

Given $U_1\in \mathcal{F}_1$, we know that there are infinitely many $j\ge g(U_1)$ such that
\begin{eqnarray*}
\#\mathcal{Q}(U_1,j,\delta(U_1))   & \ge &    \kappa  \cdot  2^{d(j-g(U_1))}\\
 \#\widetilde {\mathcal{Q}}(U_1,j,\delta(U_1))  & \le &     C_\delta  \cdot 2^{d(j-g(U_1))}  \Big( 2^{-dj\varphi(2^{-j})} \\ && + \sum_{g(U_1)\le k\le j/\delta(U_1)}    2^{-dk(\delta(U_1)-1-\psi(2^{-k}))}
\Big).
\end{eqnarray*}
The arguments used  in the first step to find an upper bound for the sum in the second inequality above also apply here. When $j$ is chosen   large enough, we can find a subset $\widetilde {\mathcal{F}}_2(U_1)$ of cardinality $\displaystyle \frac{ \kappa}{2^{d+1}} 2^{d(j-g(U_1))}$ in $\mathcal{Q}(U_1,j,\delta(U_1))\setminus \widetilde {\mathcal{Q}}(U_1,j,\delta(U_1))$ such that
\begin{itemize}
\item
 the  dyadic cubes $\widetilde V$ belonging to  $\widetilde {\mathcal{F}}_2(U_1)$  are mutually distant from at least $2^{-j}$, 
\item
each dyadic cube $\widetilde V\in  \widetilde {\mathcal{F}}_2(U_1)$ satisfies simultaneously $\widetilde V\in \mathcal{G}_j$, $\widetilde V\subset U_1$,  
$
 \mathcal{P}\big (\widetilde V,\delta(U_1)\big )  
  $ and 
  $$ \widetilde V \ \bigcap \left(  \  
\bigcup_{k=g(U_1)}^{\gamma(j)} \  \bigcup_{p\in \mathcal{T}_k} \ B(y_p, (\rho_p)^{\delta(U_1)}) \right ) = \emptyset.
$$
\end{itemize}
As in the first step, we associate with every $\widetilde V\in  \widetilde {\mathcal{F}}_2(U_1)$ a dyadic cube called its contracted descendant $V$, which enjoys the following properties:
\begin{itemize}
\item
There exists an element $x(\widetilde V) \in \widetilde V$ and a  positive real number $r(\widetilde V)$ such that $(x(\widetilde V),r(\widetilde V)) \in \mcs$, $r(\widetilde V)$ satisfies $2^{-j-1}\le r(\widetilde V)\le 2^{-j}$,  and every $x\in V$ is    approximated at a rate belonging to $[\delta(U_1) ,\delta(U_1)+\varphi(2^{-j})]$ by $(x(\widetilde V),r(\widetilde V)) $,
  \item
  if $p$ is such that $2^{-g(V)}\le r_p\le 2^{-g(U_1)}$ and $x_p\neq x_n$, then for all $x\in  V$ we have $x \notin  B\big (x_p, (r_p)^{\delta(U_1)+\theta (2^{-j})H_\ep/2}/16\big )$.
\end{itemize}

We now fix the integer $j=j({U_1})$ so that all the assumptions above are satisfied, and we set 
$$
\mathcal{F}_2(U_1)=\{V: \mbox{$V$ is the contracted descendant of one }\widetilde V\in \widetilde {\mathcal{F}}_2(U_1)\},
$$
and 
$$
\mathcal{F}_2=\{ V \in \mathcal{G}: \  \exists \,  U_1 \in\mathcal{F}_1 \mbox{ such that } V\in \mathcal{F}_2(U_1)\}.
$$
 The measure $\mu_\varepsilon$ can be extended into a Borel probability measure on the algebra  $\sigma_{2}=\sigma ( L: L\in \mathcal{F}_1\cup \mathcal{F}_2)$ by imposing
$$
\mbox{for every $U_1\in \mathcal{F}_1$, for every $V\in \mathcal{F}_2(U_1)$, }  \ \ 
\mu_\varepsilon(V)=\frac{\mu_\varepsilon(U_1)}{\# \mathcal{F}_2(U_1)} .
$$

We choose $j_1:=\min(j(U_1): U_1 \in \mathcal{F}_1)$ large enough so that for every $U_1 \in \mathcal{F}_1$, $ j_1 \ge 2g(U_1) $ and 
\begin{eqnarray}
\label{eq5}
\max \left( \frac{2^{d+1}} {\kappa} , 2^{dg(U_1)} , |U_1|^{d/\delta(U_0) -3d\ph(|U_1|)}\right) \le 2^{dj_1\varphi(2^{-j_1})}.
 \end{eqnarray}
In particular, for every $V\in\mathcal{F}_2(U_1)$ we have $\frac{2^{d+1}} {\kappa}\le |V|^{-d\varphi(|V|)}$. Moreover,   for some constant $C>0$, 
\begin{eqnarray}\label{eq6}
 \mbox{for every $V\in\mathcal{F}_2(U_1)$, }  \ \   |V|\ge C \cdot {2^{-dj(U_1) \delta(U_1) }}.
 \end{eqnarray}
 
Let us check the scaling properties of the measure $\mu_\varepsilon$ on the elements of $\sigma_2$. Let $U_1\in \mathcal{F}_1$ and $V\in \mathcal{F}_2(U_1)$. Combining  (\ref{maj1}) and the lower bound for the cardinality of $\widetilde{ \mathcal{F}}_2(U_1)$, we obtain  that
\begin{eqnarray*}
\mu_\varepsilon(V) &\le& \displaystyle \frac{2^{d+1}} { \kappa}  \cdot 2^{d(g(U_1)-j(U_1))}  \mu_\varepsilon(U_1) \\
&\leq&\displaystyle \frac{2^{d+1}}{ \kappa}  \cdot 2^{d(g(U_1)-j(U_1))} |U_1|^{d/\delta(U_0) -3d\ph(|U_1|)}.
\end{eqnarray*}
By (\ref{eq5}) and then (\ref{eq6}), we get
\begin{eqnarray*}
\mu_\varepsilon(V)   \le  2^{ -dj(U_1)} 2^{- 3 dj_1\ph(2^{-j_1})} \leq C  |V| ^{ d/ \delta(U_1) } 2^{- 3d j_1\ph(2^{-j_1})}.
\end{eqnarray*}
Using the monotonicity of $r\mapsto r^{-\ph(r)}$, which  tends to $+\infty$ when $r\to 0^+$, we see that $|V|^{-\ph(|V|)} \geq 2^{- j_1\ph(2^{-j_1})}$ when $j_1$ is large enough. We get 
\begin{equation}
\label{maj2}
\mu_\varepsilon(V) \leq |V| ^{ d/\delta(U_1)- 3d\varphi( |V|)}.
\end{equation}

As in the first step, given $U_1\in\mathcal{F}_1$, for any pair of distinct elements of $\mathcal{F}_2(U_1)$, namely $(V,V')$, we have $d(V,V')\le 2^{-j(U_1)}$.

\subsection{Induction}

Suppose that for $n\ge 2$ we have constructed $\mathcal{F}_{0},\dots,\mathcal{F}_n$, a finite sequence of sets of closed dyadic cubes, as well as a measure $\mu_\varepsilon$ on $\sigma_n=\sigma\Big (I:I\in\bigcup_{1\le m\le n}\mathcal{F}_m \Big )$ such that:

\begin{enumerate}
\item
 For every $1\le m\le n$, each element $U$ of $\mathcal{F}_m$ is included in one element of $\mathcal{F}_{m-1}$, and satisfies $ \displaystyle \frac{2^{d+1}} {\kappa}  \le  |U|^{-d\varphi(|U|)}$.

\medskip
\item For every $1\le m\le n$, if $ U\in \mathcal{F}_{m-1}$, there exists a dyadic generation $j(U)  $ such that:

$(a)$ We have
\begin{equation}
\label{eq8}
\!\! 2g(U)   \le   j (U)  \mbox{ and }    \# \{V\in\mathcal{F}_m: V\subset U\}= \displaystyle \frac{\kappa \cdot  2^{d(j(U)-g(U))}}{2^{d+1}}, 
\end{equation}
and if two distinct elements $V$ and $V'$ of $\mathcal{F}_m$ belong to $U$ then $d(V,V')\geq 2^{-j(U)}$. 

\sk 

$(b)$ for every $V \in \mathcal{F}_m$ such that $V\subset U$, there exist a cube $\wv\in \mathcal{Q}(U,j(U),\delta(U))\setminus \widetilde {\mathcal{Q}}(U,j(U),\delta(U))$ such that $V\subset \widetilde V\subset U$, as well as an element $x(\wv)\in\wv$ and a positive real number $r(\wv)$ satisfying  $(x(\wv),r(\wv))\in \mcs $ and $2^{-j(U)-1}\le r(\wv)\le 2^{-j(U)}$. Moreover,  every  element $x\in V$ is approximated at a rate belonging to $[\delta(U) , \delta(U)+\varphi(2^{-j(U)})]$ by $(x(\wv),r(\wv))$. 

$(c)$  if $p\ge 1$ satisfies $g(U)\le r_p\le g(V)$ and $x_p\neq x (\wv)$, then  no element  $x\in V$ belongs  to $ B\big (x_p,r_p^{\delta(U)+\theta(2^{-j(U)})H_\ep/2}/16\big )$.

\sk

\item
 If $1\le m\le n$ and $U\in\mathcal{F}_{m-1}$, then for $V\in\mathcal{F}_{m}$ such that $V\subset U$ we have:
$$
\mu_\varepsilon( V)=\frac{\mu_\varepsilon(U)}{\# \{V' \in \mathcal{F}_{m}: \  V' \subset U\}}.
$$

\item
 For all $1\le m\le n$,  $U\in \mathcal{F}_{m-1}$ then for  $V\in\mathcal{F}_m$ such that $V\subset U$ we have
$$
\mu_\varepsilon(V)\le   |V|^{d/\delta(U) -3d\varphi (|V|)}.
$$
\end{enumerate}

\ms

Parts 1. to 4.  of the induction are easily checked for the first  generations $\mathcal{F}_1$ and $\mathcal{F}_2$.

The technique we use to build the   generation $\mathcal{F}_{n+1}$  is the same as for the first iteration. We briefly indicate the steps to follow.

\ms

For each $  U_n\in\mathcal{F}_n$, we know that  there are infinitely many integers $j\ge g(U_n)$ such that
\begin{eqnarray*}
\#\mathcal{Q}(U_n,j,\delta({U_n}))    & \ge &    \kappa \cdot  2^{d(j-g(U_n))}\\
 \#\widetilde {\mathcal{Q}}(U_n,j,\delta({U_n}))  & \le &    C_d \cdot 2^{d(j-g(U_n))}  \Big( 2^{-dj\varphi(2^{-j})} \\ && + \sum_{g(U_n)\le k\le j/\delta(U_n)}    2^{-dk(\delta(U_n)-1-\psi(2^{-k}))}
\Big).
\end{eqnarray*}
If the integer  $j=j(U_n)$ is chosen large enough, there is a set $\widetilde {\mathcal{F}}_{n+1}(U_n)$ of cardinality $\displaystyle \frac{\kappa}{2^{d+1}} \cdot 2^{d(j-g(U_n))}$ included in  $\mathcal{Q}(U_n,j,\delta({U_n}))\setminus \widetilde {\mathcal{Q}}(U_n,j,\delta({U_n}))$ such that
\begin{itemize}
\item
 the  dyadic cubes $\widetilde V$ belonging to  $\widetilde {\mathcal{F}}_{n+1}(U_n)$  are mutually distant from at least $2^{-j(U_n)}$, 
\item
each $\widetilde V\in  \widetilde {\mathcal{F}}_{n+1}(U_n)$ satisfies simultaneously $\widetilde V\in \mathcal{G}_j$, $\widetilde V\subset U_n$,   
$ \mathcal{P}\big (\widetilde V,\delta({U_n}) \big )  $ and 
  \begin{eqnarray*}  \widetilde V \ \bigcap \left(  \ 
\bigcup_{k=g(U_n)}^{\gamma (j)}  \ \bigcup_{p\in \mathcal{T}_k} \ B(y_p, (\rho_p)^{\delta({U_n})}) \right ) = \emptyset.
\end{eqnarray*}
\end{itemize}
We can associate with each  $\widetilde V\in  \widetilde {\mathcal{F}}_{n+1}(U_n)$   a contracted descendant $V$, which is a dyadic cube enjoying the properties:  
\begin{itemize}
\item
By condition $\mathcal{C}_2$, there is $x(\wv)\in\wv$ and a positive real number $r(\wv)$ satisfying  $(x(\wv),r(\wv))\in \mcs $ and $2^{-j(U_n)-1}\le r(\wv)\le 2^{-j(U_n)}$. Moreover,  every  element $x\in V$ is approximated at a rate belonging to $[\delta(U_n) , \delta(U_n)+\varphi(2^{-j(U_n)})]$ by $(x(\wv),r(\wv))$.     
  \item
  if $p\ge 1$ is such that $2^{-g(V)}\le r_p\le 2^{-g(U_n)}$ and $x_p\neq x(\wv)$, then for all $x\in  V$ we have $x \notin  B\big (x_p, (r_p)^{\delta({U_n})+\theta (2^{-j(U_n)})H_\ep/2}/16\big )$.
\end{itemize}
Then we set 
$$
\mathcal{F}_{n+1}(U_n)=\{V: \mbox{$V$ is the contracted descendant of some  }\widetilde V\in \widetilde {\mathcal{F}}_{n+1}(U_n)\},
$$
and 
$$
\mathcal{F}_{n+1}=\{ V \in \mathcal{G}: \  \exists \,  U_n \in\mathcal{F}_{n+1} \mbox{ such that } V\in \mathcal{F}_{n+1}(U_n)\}.
$$

The measure $\mu_\varepsilon$ can be extended into a Borel probability measure on the algebra  $\sigma_{n+1}=\sigma ( L: L\in \bigcup_{p=0}^{n+1} \mathcal{F}_p)$ by the following formula:
$$
\mbox{for every $U \in \mathcal{F}_{n+1}$, for every $V\in \mathcal{F}_{n+1}(U)$, }  \ \ 
\mu_\varepsilon(V)=\frac{\mu_\varepsilon(U)}{\# \mathcal{F}_{n+1}(U)} .
$$
In addition, requiring that   $j_n:=\min(j(U_n): U_n \in \mathcal{F}_n)$ is  large enough so that for all $U \in \mathcal{F}_n$ and $ T\in\mathcal{F}_{n-1}$ such that $U\subset T$,  we obtain that  $j(U) \ge 2g(U)$ and   
\begin{eqnarray*}
\max \left(  \frac{2^{d+1}} {\kappa}  , 2^{dg(U)} , |U|^{d/\delta({T}) - 3d\ph(|U|)}\right) \le 2^{dj_n\varphi(2^{-j_n})}.
 \end{eqnarray*}
This  ensures that (\ref{eq8})   holds with $p=n+1$. Finally the same lines of computations as in the second step of the construction yield the part 4.  of the induction, i.e. the scaling behavior of the measure $\mu_\varepsilon$ on the dyadic cubes of the $(n+1)$th generation of the Cantor set.

\ms

Iterating the previous construction,  the Kolmogorov extension theorem yield a measure $\mu_\varepsilon$ on the algebra $\sigma\Big (V: \ V\in\bigcup_{n\ge 1}\mathcal{F}_n \Big )$ such that all the properties 1. to 4.  hold true for all $n\ge 1$. By construction, the measure $\mu_\varepsilon$ is carried by the Cantor set $$\mathscr{K}_\varepsilon=\bigcap_{n\ge 0} \ \bigcup_{V\in\mathcal{F}_n} \ V.$$

\subsection{Scaling properties of $\mu_\varepsilon$}\label{scalingprop}

Let $\delta_\varepsilon=\sup_{U\in \bigcup_{n\ge 0}\mathcal{F}_n} \delta (U)$. By (\ref{eq00}), $\de_\ep \leq H_\ep:=f(x_\ep)+ \ep(2 \frac{d}{h^2} +1)$. 

We are going to show that there exists $C'>0$ such that for every open cube $B\subset [0,1]$,
\begin{equation}\label {dim}
\mu_\varepsilon (B)\le C'  |B|^{d/\delta_\varepsilon}| B|^{-4d\varphi (|B|)}.
\end{equation}
If (\ref{dim}) holds true, then Lemma \ref{mdp}, known as the {\em mass distribution principle} \cite{F2},  allows
to bound by below  the Hausdorff dimension of the support of $\mu_\varepsilon$.
\begin{lemma}
\label{mdp}
Let $F$ be a Borel set in $\R^d$, and $\mu$ be  a Borel probability measure on
$F$. Suppose that, for some $\eta>0$, there are 
$\alpha>0$ and a gauge function $\zeta$ such that 
$\liminf_{x\to 0^+} \frac{\zeta(x)}{x^\alpha} >0
$ and for every set $U$ with a diameter less than $ \eta$,   $\mu(U) \leq C
\zeta(|U|).$ 

Then $ \mathcal {H}^{\zeta }(F) \geq \mu(F)/C$ and ${\dim_{\mathcal{H}}}
F\geq \alpha$.
\end{lemma}

 Let $B$ be an open subcube of $[0,1]^d$  intersecting $\mathscr{K}_\varepsilon$. Let $n_0$ be the smallest integer such that $B$ intersects at least two elements of $\mathcal{F}_{n_0}$. By construction, the elements $V$ of $\mathcal{F}_{n_0}$  intersecting $B$ are all contained in the same element $U$ of $\mathcal{F}_{n_0-1}$, and   $\mu_\varepsilon (B)\le \mu_\varepsilon (U)$.

\ms

Suppose first that $|B|\ge |U|$. Part  4. of the induction yields  
$$
\mu_\varepsilon (B)\le \mu_\varepsilon (U)\le  |U|^{d/\delta_\varepsilon -3d \varphi (|U|)}\le   |B|^{d/\delta_\varepsilon -3d\varphi (|B|)}$$
when $|B|$ is small enough. Once again the monotonicity of $r\mapsto r^{-\ph(r)}$ has been used.

\ms

Suppose now that $|B|<| U|$. Applying Part 4. of the induction, we find 
$$
\mu_\varepsilon (B)\le \mu_\varepsilon(U)\frac{\#\{V\in\mathcal{F}_{n_0}: V \subset U,\ V\cap B\neq\emptyset\}}{\# \{V\in \mathcal{F}_{n_0}:\ V\subset U\}}.
$$

Let us use Part 2. of the induction to bound by above  $\#\{V\in\mathcal{F}_{n_0}: V \subset U,\ V\cap B\neq\emptyset\}$.  There exists an integer $j(U)$ such that the elements of $\mathcal{F}_{n_0}$ that intersect $B$ are distant from one another by at least $2^{-j(U)}$ and have   diameter less than $2^{-j(U)}$. Consequently, due to Lemma \ref{intersection}.1, there are at most $C_d|B|^d2^{dj(U)}$ of them.

 In addition, we know that 
$$\# \{ V \in \mathcal{F}_{n_0},\ V\subset U\}\ge \displaystyle \frac{\kappa}{2^{d+1}}  \cdot 2^{-dg(U)} 2^{dj(U)} \ =   \frac{\kappa}{2^{d+1}}  \cdot  |U|^d 2^{dj(U)}.$$
 This yields thanks to (\ref{eq8})
\begin{eqnarray*}
\mu_\varepsilon (B)&\le & \mu_\varepsilon(U)\frac{2|B|^d2^{dj(U)}}{ \frac{\kappa}{2^{d+1}} | U|^d 2^{dj(U)}} \leq  2  \cdot  \mu_\varepsilon(U ) \frac{|B|^d}{|U|^d}|U|^{-\varphi(|U|)}.
\end{eqnarray*}
Using the scaling behavior of $\mu_\varepsilon$ on the elements of $\mathcal{F}_{n_0}$, we get 
\begin{eqnarray*}
\mu_\varepsilon (B)&\le & 2 \cdot  |U|^{d/\delta_\varepsilon -3d\varphi (|U|)} \frac{|B|^d}{|U|^d} |U |^{-d\varphi(|U|)}\\
&\leq & 2 \cdot   |B|^{d/\delta_\varepsilon -4d \varphi (|B|)} \Big (\frac{|B|}{|U|}\Big )^{d-d/\delta_\varepsilon }\frac{|B|^{4d\varphi (|B|)}}{|U|^{4d\varphi(|U|)}}\\
&\leq & 2  \cdot |B|^{d/\delta_\varepsilon-4d\varphi (|B|)},
\end{eqnarray*}
the last line following from the observation that 
  $\displaystyle \Big (\frac{|B|}{|U|}\Big )^{d-d/\delta_\varepsilon}\frac{|B|^{4\varphi (|B|)}}{| U|^{4\varphi(|U|)}}$ is bounded by above by  1  due to the monotonicity property of $r^{\varphi(r)}$ and the fact that $|B|<|U|$.

\ms

By the mass distribution principle, the Hausdorff dimension of $\mu_\varepsilon$ (and thus the Hausdorff dimension of $\mathscr{K}_\varepsilon $) is larger than  $\displaystyle \frac{d}{\delta_\varepsilon}$, which by (\ref{eq00}) is greater than
\begin{eqnarray}
\nonumber  \frac{d}{\delta_\ep}  & \geq & \frac{d}{H_\ep}\geq \frac{\de}{f(y_\ep)+\ep (2\frac{d}{h^2}+1)} \geq \frac{d}{f(y_\ep)} \cdot  \frac{1}{1+\ep (2\frac{d}{h^2}+1)/{f(y_\ep)}}\\
\label{defheps}
& \geq & (h-\ep)  \cdot \frac{1}{1+\ep (2\frac{d}{h^2}+1)/{f(y_\ep)}} := h_\ep.
\end{eqnarray}
It is obvious that $h_\ep $ increases toward $ h$ when $\ep$ goes to zero, hence the result.

\subsection{Relation with $  \ws(\Omega \setminus E,f) $}  

  Let us prove that $\mathscr{K}_\varepsilon\setminus E\subset \ws(\Omega \setminus  E,f)$. 
  
  \sk
  
Let $x\in \mathscr{K}_\ep\setminus E$ and for $n\ge 1$ denote by $U_n(x)$ the unique element of $\mathcal{F}_n$ that contains $x$. Using parts 2. and 3. of the induction, we have  $\delta_x=\limsup_{n\to\infty}\delta(U_n(x))$.   
    
Recall   that the function $f$ is continuous at $x$. Using  formula (\ref{defdeltau}), one observes that  $f(y_{U_n(x)})$ converges to $f(x)$ (since $y_{U_n(x)}$ is any point of $U_n(x)$). This implies that  $\delta(U_n(x))$   converges to  $\displaystyle  {f(x)}$ when $n$ tends to infinity.

\sk

Finally,  $\delta_x=  \lim_{n\to\infty}\delta(U_n(x))  =f(x)$.

\subsection{Density of $\ws(\Omega \setminus E,f) $ in $\Omega$ when ${\dim_{\mathcal{H}}}\, E< \frac{d}{\sup \{f(x): x\in \Omega\setminus E\}}$}
 
 Using what precedes, we are able to construct a Cantor set $\mathcal{K}_\ep$ in order to approximate the Hausdorff dimension of $\ws(\Omega  \setminus E,f) $. But our construction may be achieved in a neighborhood $U_y$ of any point $y\in  \Omega\setminus E$ such that ${\dim_{\mathcal{H}}}\, E\cap U_y < \frac{d}{\inf \{f(x): x\in U_y \setminus E\}}$. Consequently, if ${\dim_{\mathcal{H}}}\, E< \frac{d}{\sup \{f(x): x\in \Omega\setminus E\}}$ then we get the conclusion, since $\Omega\setminus E$ is dense in $\Omega$. 

\section{Upper bounds for the dimensions}
\label{secmaj}

 We suppose that the assumptions of Theorem~\ref{main2} are fulfilled. As in the previous section, we  set  $\de=\inf \{f(x): x\in \Omega\setminus E\}$ and $h= d/ \de$.
  
 \ms

By (\ref{cond1}),  we know that $\delta_x \ge 1$ for all $x\in\Omega$. The set $\ws(\Omega \setminus E,f)$ contains only elements $x\in \zu^d$ satisfying $f(x)\ge \de$, which implies that $\delta(x)\ge \de$. Hence, for every $\ep>0$, 
 $ \ws(\Omega \setminus E,f) \subset \mathcal{L}_{\de-\ep}$, where we recall that $ \mathcal{L}_{\de-\ep}$ is given by (\ref{defjarnik0}):
$$ \mathcal{L}_{\de-\ep}= \bigcap_{N\geq 1} \, \bigcup_{n\geq N} B(x_n,r_n^{\de-\ep})=\bigcap_{N\geq 1} \, \bigcup_{n\geq N} B(y_n,\rho_n^{\de-\ep}).$$

It is known (see   \cite{BSUbiquity}) that if the system $\mcs$ satisfies  $\mathcal{C}_1$, then ${\dim_{\mathcal{H}}}\, \mathcal{L}_{\delta }\le d/\de $ for all $\de\ge 1$.  Let us prove it briefly   for completeness.

 Let $s> d/\delta$.  For any integer $N\geq 1$, a covering of the limsup set $\mathcal{L}_{\de}$ is provided  by the union of sets $\bigcup_{n\geq N} B(y_n,\rho_n^{\de})$. Let $\eta>0$, and choose $N$ large enough so that $2\rho_n\leq \eta$ for $n\geq N$. Recalling the definition of the generalized Hausdorff measure associated with the gauge function $\zeta_s(x)=x^s$, we see that 
$$
  \mathcal{H}_\eta^{\zeta_s}(\mathcal{L}_{\de})  \leq \sum_{n\geq N} |B(y_n,\rho_n^{\de}) |^s  \leq \sum_{j=J}^{+\infty} \ \sum_{p\in \mathcal{T}_j}  |B(y_p,\rho_p^{\de})|^s \leq  \sum_{j=J}^{+\infty} \ \sum_{p\in \mathcal{T}_j}  2^{-js\de} ,
$$
where $J$ is the unique integer such that $y_N \in \mathcal{T}_J$. Using Condition $\mathcal{C}_1$, we see that 
 \begin{eqnarray*}
  \mathcal{H}_\eta^{\zeta_s}(\mathcal{L}_{\de})  & \leq &  \sum_{j=J}^{+\infty} \ N_j \cdot 2^{dj-js\de} \leq \sum_{j=J}^{+\infty} \ N_j \cdot 2^{j(d-s\de)}.
  \end{eqnarray*}
This series converges, since $\log (N_j)=o(j)$ and $d-s\de<0 $ by construction. Consequently, the  $s$-Hausdorff measure of $\mathcal{L}_{\de }$ is finite for any $s > d/\delta$. This demonstrates that   $\dim_{\mathcal{H}} \,\mathcal{L}_{\de} \leq  d/\de$.

\ms

The above argument applies to $\mathcal{L}_{\delta-\ep}$ when $\de-\ep>1$, and thus $${\dim_{\mathcal{H}}}  \ws(\Omega \setminus E,f) \leq \inf_{\ep>0} {\dim_{\mathcal{H}}}\, \mathcal{L}_{\delta-\ep} \leq  \inf_{\ep>0}  d/(\delta-\ep) = d/\delta =h.$$
This yields the conclusion.

\section{Examples of suitable systems $
\big((x_n,r_n)
\big)_{n\ge 1}$}
\label{secexa}

\subsection{Approximation by $b$-adic numbers}  

We prove that the dyadic system satisfies  $\mathcal{C}_1$ and $\mathcal{C}_2$. The case of the  $b$-adic system (whose definition is clear)  is similar.

\sk

Define the system $\mathcal{D}=  \big(({\bf k} \cdot 2^{-j}, 2^{-j})
\big)_{j\geq 1, {\bf k}\in\{0,1,...,2^{-j}-1\}^d}$, and consider the approximation rate of any $x\in \zu^d$ by  $\mathcal{D}$
$$\de_x= \sup\{\de\geq 1 : |x-{\bf k}\cdot 2^{-j}|\leq 2^{-j\de} \mbox{ for an infinite number of $(j,{\bf k})$}\}.$$
We rather consider the system $\mathcal{D}'= \big(({\bf k}\cdot 2^{-j}, \displaystyle\frac{2^{-j}}{32})\big)_{j\geq 1, {\bf k}\in\{0, ...,2^{-j}-1\}^d}$ and the associated approximation rate
$$\de'_x= \sup\{\de\geq 1 : |x-{\bf k}\cdot 2^{-j}|\leq ( \displaystyle\frac{2^{-j}}{32})^\de \mbox{ for an infinite number of $(j,{\bf k})$}\}.$$
Of course, $\de_x=\de'_x$ for every $x\in\zu^d$, but the   constant
$32$ is necessary for our condition $\mathcal{C}_2$ to hold.
 
 \sk
 
 The irreducible subsystem of $\mathcal{D}'$  consists in the couples $({\bf k} \cdot 2^{-j},\displaystyle\frac{2^{-j}}{32})$ for which at least one coordinate of ${\bf k}$ is odd. Therefore, it is obvious that the weak redundancy condition $\mathcal{C}_1$ is satisfied,  the corresponding sequence $(N_j)_{j\geq 1}$ being constant equal to 1, so that $\gamma(j)=j$ for every $j\geq 1$. 
 
 
 \sk
 
 To check  $\mathcal{C}_2$, let $\de>1$, and  consider any $\ph\in\Phi$. Let  ${\bf k} \cdot 2^{-j}$  be a dyadic element of $\zu^d$ such that   ${\bf k}$ has at least one odd coordinate. We call $V$ the dyadic cube $\prod_{i=1}^d [k_i\cdot 2^{-j}, (k_i+1)\cdot 2^{-j})$  of generation $g(V)=j$.  Given a dyadic generation $j$, the number of such dyadic irreducible cubes is greater than $2^{dj-1}$.  Then the property $\mathcal{P}(V,\de )$ holds   without any further condition. Indeed, we only have to check that for every $j'\in \{\gamma(j)j,...,(j+1) \de +4\}$, for every ${\bf k}'$ (with at least one odd coordinate), ${\bf k}'2^{-j'} \notin B({\bf k} \cdot 2^{-j},(\displaystyle\frac{2^{-j}}{32})^{\de })$. This is obvious, since by  the structure of the dyadic tree we get when $j\leq j'\leq (j+1) \de +4$
 \begin{eqnarray*}
 |{\bf k}\cdot 2^{-j} -{\bf k}'\cdot 2^{-j'} | & \geq & 2^{-j'}  \geq 2^{-(j+1) \de  -4} \geq \frac{2^{-j  \de  }}{ 16\cdot 2^{\de }}\geq  (\displaystyle\frac{2^{-j}}{32})^{\de  },
 \end{eqnarray*}
  and when $\gamma(j)\leq j'<j$
  \begin{eqnarray*}
 |{\bf k}\cdot 2^{-j} -{\bf k}'\cdot 2^{-j'} | & \geq & 2^{- j}   \geq  (\displaystyle\frac{2^{-j}}{32})^{\de  }.
 \end{eqnarray*}
 
 Thus the system $\mathcal{D}$ satisfies $\mathcal{C}_2$ with a  function $\kappa$ constant equal to 1/2 (it holds for all the "irreducible" sub-cubes of $\zu^d$).

\subsection{Diophantine approximation by  rational numbers in $\R$}  

Consider the system 
$$\mathcal{R} = \big ((p/q, 1/q^{2}) \big)_{q\geq 1, \, 0\le p\le q-1}.$$
 It follows from Dirichlet's argument   that $ \mathcal{L}_1(\mathcal{R}) =[0,1] $. The irreducible sub-system of $\mathcal{R}$ consists in the elements  of $ \mathcal{R}$ such that $p \land q=1$.

\ms

We are going to check that $\mathcal{R}$ satisfies $\mathcal{C}_1$ and $\mathcal{C}_2$.

\ms

Let $j\ge 1$ be an integer, and let $(p/q,1/q^{2})\in\mathcal{R} $ be   such that $q^{2} \in (2^{j},2^{j+1}]$. We shall prove that  $B({p}/q,1/q^{2})$ may contain only a bounded number of rational numbers $p'/q'$ satisfying  $( {p}'/q',1/(q')^{2})\in\mathcal{R} $ and  $(q')^{2} \in (2^{j},2^{j+1}]$.   This  implies  $\mathcal{C}_1$.

If $p /q\neq p' /q'$, then one has  necessarily that $|p /q-p' /q'|= |pq'-p'q|/(qq')\ge 1/(qq') \geq 2^{-j-1}$, since $q$ and $q'$ belong to $[2^{j/2}, 2^{(j+1)/2})$.  Since the diameter of $B({p}/q,1/q^{2})$ is at most $2^{-j+1}$,   there are at most $4$ distinct irreducible rational numbers  $p'/q'$  belonging to $B(p/q, 1/q^2)$. Hence   $\mathcal{R}$ satisfies $\mathcal{C}_1$, with a sequence  $(N_j)_{j\ge1 }$ constant  equal to 4. 
 
\ms

In order to prove $\mathcal{C}_2$, we consider   $V:=[K\cdot 2^{-J}, (K+1)\cdot 2^{-J}) $  a dyadic interval in $\zu$ of generation $J$, a real number $\de>1$ and any function $\ph\in\Phi$. We demonstrate that $\mathcal{P}(V,\de)$ holds without any restriction on $V$, $\de$ and $\ph$. Obviously $V$  contains a rational number  ${p}/q$ satisfying  ${q}^{2}\in (2^{J},2^{J+1}]$ ($p/q$ is not necessarily irreducible). Assume that a rational number  $p'/q' \neq p/q$ belongs to  $B(p/q, 1/q^{2\de})$ with $ \log_2 ((q')^2)  \in [{\gamma(J)}, \cdots, {\de(J+1) +4}]$. This implies that $q'  \leq 2^{({\de( J+1) +4})/2}   \leq q^{\de/2} 2 ^{\de/2+2 } $.  Combining the information, we have
$$1/q^{2\de} \geq |p/q-p'/q'| \geq 1/(qq') \geq 2 ^{-\de/2-2 }/q ^{1+\de/2}.$$
This last inequalities can not hold as soon as $\de >1$ (provided that $q$ is large enough). Consequently, $\mathcal{P}(V,\de)$  holds, and $\mathcal{R}$ satisfies  $\mathcal{C}_2$ with a  function $\kappa$ constant equal to 1.
  
\subsection{Inhomogeneous Diophantine approximation}  

Let $\alpha$ be an irrational number in $\zu$.
Consider the system
$$\mathcal{I}= \Big(\big(\{n\alpha\}, \frac{1}{n}\big)\Big)_{n\geq 1},$$
where $\{x\}$ stands for the fractional part of the real number $x$.

It is proved in \cite{BSUbiquity} (Proposition 6.1) that $\mathcal{I}$ satisfies $\mathcal{C}_1$ if and only if the approximation rate of $\alpha$  by the rational system $\mathcal{R}$ equals 2.

We prove that $\mathcal{I}$ satisfies $\mathcal{C}_2$, when the approximation rate of $\alpha$ by  the  rational system $\mathcal{R}$ is 2. When this holds, for every $\ep>0$, there is an integer $q_\ep$ such that
\begin{equation}
\label{defqep}
\mbox{for every $q\geq q_\ep$, \  for every integer $p$, } \  |\alpha -p/q| \geq 1/q^{1+\ep}.
\end{equation}
We focus now on $\mathcal{C}_2$. For this, let us recall the {\em three distance theorem} \cite{SOS,SLATER,BUGEAUD3}: the real numbers $\{\alpha\}$, $\{2\alpha\}$, $\{3\alpha\}$, ..., $\{ N\alpha\}$ divide the interval $\zu$ into $N+1$ intervals whose lengths take at most three values $d_1(N)$, $d_2(N)$ and $d_3(N)$, satisfying 
$$d_1(N)<d_2(N)<d_3(N) \leq \frac{3}{N+1}.$$
Let $J\geq 1$. As for the rational system, in order to prove $\mathcal{C}_2$, we consider   $V:=[K\cdot 2^{-J}, (K+1)\cdot 2^{-J}) $  a dyadic interval in $\zu$ of generation $J$, a real number $\de>1$ and any function $\ph\in\Phi$.  We demonstrate that $\mathcal{P}(V,\de)$ holds without any restriction on $V$, $\de$ and $\ph$ for a sufficiently large number of dyadic intervals $V$. 

 Apply the three distance theorem to $\{\alpha\}$, $\{2\alpha\}$, $\{3\alpha\}$, ..., $\{ 2^{J}\alpha\}$. The $2^{J}+1 $ corresponding intervals of $\zu$ have length less than $3/(2^{J}+1)$. By a translation argument, the points $\{( 2^{J}+1)\alpha\}$, $\{( 2^{J}+2)\alpha\}$, $\{( 2^{J}+3)\alpha\}$, ..., $\{ 2^{J+1}\alpha\}$ divide the interval $\zu$ into $2^J+1$ intervals whose lengths are also less than $3/(2^{J}+1)$. This means that among the dyadic intervals of generation $J$, there are no three consecutive dyadic intervals $U$ which do not contain one of the points $\{n\alpha\}$, for $n$ ranging over $\{2^J+1,2^J+2,..., 2^{J+1}\}$.  
 
 Let us consider one such interval $V:=[K\cdot 2^{-J}, (K+1)\cdot 2^{-J}) $, which contains $\{n\alpha\}$ for some  $n$ belonging to $\{2^J+1,2^J+2,..., 2^{J+1}\}$. Assume that another point $\{n'\alpha\}$ belongs to $B(\{n\alpha\},1/n ^\delta)$ with $\log_2 n' \in [\gamma(J),\cdots,  [\de (J+1)] +4]$. This means that
$ |\{n'\alpha\}- \{n\alpha\}| \leq \frac{1}{n^\de}$.  By definition there are integers $p$ and $p'$ satisfying
$n\alpha = p+ \{n'\alpha\}$ and $n'\alpha = p'+ \{n\alpha\}$, hence
$$ |(n'\alpha+p') - (n\alpha - p)|  = | (n-n')\alpha- (p' - p)| \leq \frac{1}{n^\de},$$
 or equivalently
 $$  \Big|  \alpha-  \frac{p' - p}{n'-n} \Big| \leq \frac{1}{|n'-n|\cdot n^\de} \leq  \frac{C}{|n'-n| ^{\de +1}},$$
the last inequality following from the fact that $|n'-n| \leq  2\de \cdot n$. This contradicts (\ref{defqep}).
 Consequently, $\mathcal{P}(V,\de)$  holds.
 
 Finally,  $\mathcal{I}$ satisfies  $\mathcal{C}_2$ with a function $\kappa$ constant equal to  $1/3$.

\subsection{Poisson point process}

Let $\mathcal P$ be a Poisson point process with intensity 
\begin{equation}
\label{defpoisson}
\displaystyle \Lambda=\mathbf{1}_{[0,1]\times (0,1)}(x,y)\cdot \ell(dx)\otimes \frac{\ell(dy)}{y^2},
\end{equation}
where $\ell$ stands for the Lebesgue  measure on $(0,1)$.   We rewrite $\mathcal{P}$ as $\mathcal{P}= \big((x_n,r_n)\big)_{n\ge 1}$, where $(r_n )_{n\geq 1}$ is a positive decreasing sequence converging to zero when $n$ tends to infinity. 

\ms

With probability one, such a system satisfies $\mathcal{C}_1$, see for instance Proposition 6.2 in   \cite{BSUbiquity}.

\medskip

We now deal with $\mathcal{C}_2$. We only need to find  a function $\varphi\in\Phi$ and a continuous function $\kappa:(1,+\infty) \to \R_+^*$ such that for every $\delta>1$,  with probability 1, for every $U\in\mathcal{G}_*$, there are infinitely many integers $j\ge g(U)$ satisfying  $\#\mathcal{Q}(U,j,\delta)\ge \kappa(\de) \cdot 2^{j-g(U)}$. Then, for any countable and dense subset $\Delta$ of $(1,\infty)$, with probability 1, for every $\delta\in\Delta$, for every $U\in\mathcal{G}_*$, there are infinitely many integers $j\ge g(U)$ satisfying  $\#\mathcal{Q}(U,j,\delta)\ge \kappa(\de) \cdot 2^{j-g(U)}$.

In fact, any $\varphi\in\Phi$ is suitable. 
\medskip

Let $\varphi \in \Phi$ and $\delta>1$. For $U\in\mathcal{G}_*$ and $V \subset U$ such that $V\in  \bigcup_{j >  g(U)}\mathcal{G}_j$, let us introduce the event
$$
\mathcal{A}(U,V,\delta)= \left .\begin{cases}   \exists\ n\in \mathcal{T}_{g(V)} \mbox{ such that } \ x_n\in V \mbox{ and }\\    
B(x_n,(r_n)^{\delta }) \bigcap  \ \left(\bigcup_{\gamma({g(V)}) \le k \le h(V)} \mathcal{T}_j\right) =\{x_n\}
\end{cases} \!\!\!\!\!\right \}
$$
where $h(V)= \big [\delta (g(V)+1)  \big ]+4$.  Recall that $n\in \mathcal{T}_{g(V)} $ means that $2^{-g(V)-1} < r_n \leq 2^{-g(V)}$.
Note that by construction, we have the inclusion $\mathcal{A}(U,V,\delta) \subset \{\mathcal{P}(V,\delta )\text{ holds}\}$.

For every $j\ge 1$, let $\widetilde {\mathcal{G}}_j=\big\{[2k\cdot 2^{-j},(2k+1)\cdot 2^{-j}]: 0\le k\le 2^{j}-1\big\}$.  The restrictions of the Poisson point process to the strips $V\times (0,1)$, where $V$ describes $\widetilde {\mathcal{G}}_j$, are independent. Consequently, the events $\mathcal{A}(U,V,\delta)$, when $V \in \widetilde {\mathcal{G}}_j$  and $V\subset U$, are independent (we must separate the intervals in $\widetilde {\mathcal{G}}_j$   because if $V\in \mathcal{G}_j$, $x_n\in V$ and $r_n\le 2^{-j}$, then $B(x_n,(r_n)^{\delta })$ may overlap  with the neighbors of $V$).  

We denote by $X(U,V,\delta)$ the random variable  $\mathbf{1}_{\mathcal{A}(U,V,\delta)}$. For a given generation $j>g(U)$, the random variables $(X(U,V,\delta))_{V\in \widetilde{\mathcal{G}}_j}$ are i.i.d Bernoulli variables, whose common  parameter is denoted by $p_{j}(\de)$. We have the following Lemma.
\begin{lemma}\label{bbb}
There exists  a continuous function $\kappa_1:(1,+\infty) \to (0,1)$   such that for every $j\geq 1$, $p_j(\de)\ge\kappa_1(\de)$.
\end{lemma}

Let us assume Lemma \ref{bbb} for the moment. By definition we have
$$
\#\mathcal{Q}(U,j,\delta)\ge \sum_{V\in \widetilde{\mathcal{G}}_j:\ V\subset U} X(U,V,\delta).
$$
The right hand side in the last inequality is a binomial variable of parameters $(2^{j-g(U)},p_j(\de))$, with $p_j(\de)\ge \kappa_1(\de)>0$. Consequently, there exists a constant $\kappa(\de)>0$ satisfying 
\begin{equation}
\label{eqfin}
\mathbb{P}\Big (\sum_{V\in \mathcal{G}_j,\ V\subset U} X(U,V,\delta) \ge \kappa(\de) \cdot 2^{j-g(U)}\Big )\ge 1/2
\end{equation}
provided that  $j$ large enough. The continuity of $\kappa$ with respect to the parameter $\de>1$ follows from the continuity of $\kappa_1$.

Let $(j_n)_{n\ge 1}$ be the sequence defined inductively by $j_1=g(U)+1$ and $j_{n+1}=(j_n+1)\delta +5$. We notice that the events $E_n$ defined for $n\geq 1$ by 
$$E_n = \{\#\mathcal{Q}(U, j_n,\delta)\ge \kappa(\de) \cdot  2^{j_n-g(U)}\}$$
 are independent. Moreover,   (\ref{eqfin}) implies that $\sum_{n\ge 1}\mathbb{P}(E_n)=+\infty$. The Borel-Cantelli Lemma yields that, with probability 1, there is an infinite number of generations $j_n$ satisfying  $\#\mathcal{Q}(U,  j_n,\delta)\ge \kappa(\de) \cdot  2^{j-g(U)}$. This holds true for every $U\in \mathcal{G}_*$ almost surely, hence almost surely  for every $U\in \mathcal{G}_*$. Condition $\mathcal{C}_2$ is proved.

\medskip 

We prove    Lemma~\ref{bbb}. 
For every $V\in \mathcal{G}_*$, let us introduce the sets 
$$
S_V  =   V\times [2^{-(g(V)+1)},2^{- g(V)}]\ \mbox{ and } \ 
\widetilde S_V  =  V\times [2^{-h(V)},2^{-\gamma (g(V))}].
$$
 We denote by $N_V$ and $\widetilde N_V$ respectively the cardinality of $\mathcal{P}\cap S_V$ and $\mathcal{P} \cap (\widetilde S_V \setminus S_V)$. These random variables $N_V$ and $\widetilde N_V$  are independent, and  we set  $l_V=\Lambda(S_V)$ and $\widetilde l_V=\Lambda(\widetilde S_V)$ ($\Lambda$ is the intensity of the Poisson point process (\ref{defpoisson})). Due to the form of the intensity $\Lambda$,  $N_V$ and $\widetilde N_V$ are  Poisson random variables of parameter $l_V=1$ and $\widetilde l_V  = 2^{-g(V)}\big (2^{h(V)} -2^{g(V)+1}+2^{g(V)}-2^{\gamma (g(V))}\big )$ respectively. Observe  that $\widetilde l_V  \le  2^{h(V)-g(V)}$ since by definition $\gamma(g(V))\le g(V)$.
 
  We also consider   two sequences of random variables in $\R^2$  $(\xi_p=(X_p,Y_p))_{p\ge 1}$ and $(\widetilde\xi_q=(\widetilde X_q,\widetilde Y_q))_{q\ge 1}$ such that 
  \begin{eqnarray*}
 \mathcal{P}\cap S_V & = & \{\xi_p:\ 1\le p \le N_V\}\\
\mathcal{P}\cap (\widetilde S_V\setminus S_V) & = & \{\widetilde \xi_q:\ 1\le q\le \widetilde N_V\}.
\end{eqnarray*}
The event $\mathcal{A}(U,V,\delta)$ contains the event $\widetilde {\mathcal{A}}(U,V,\delta)$ defined as 
$$
\left\{
 \,  N_V=1  \mbox{ and }  
 \,  B\big (X_1, Y_1^{\delta }\big)  \, \bigcap \, \big \{ \widetilde X_q:   1\le q\le \widetilde N_V\big\}=\{X_1\}   
 \right\},
$$
where  $\xi_1=(X_1, Y_1)$. The difference between $\mathcal{A}(U,V,\delta)$ and $\widetilde {\mathcal{A}}(U,V,\delta)$ is that the latter one imposes that there is one and only one Poisson point in $S_V$.
 We have
\begin{eqnarray*}
&&\mathbb{P}(\widetilde {\mathcal{A}}(U,V,\delta))\\
 & = & \mathbb{P}\Big( \Big\{B\big (X_1, Y_1^{\delta }\big)  \bigcap  \big \{ \widetilde X_q:   1\le q\le \widetilde N_V\big\} =  \emptyset \, \Big| \, \big\{N_V=1\big\}\Big\} \Big) \\
 &&  \times   \ \mathbb{P}(\{N_V=1\})\\
& = & \mathbb{P}\Big(  \Big\{\forall \, 1\le q\le \widetilde N_V, \  \widetilde X_q \not\in B\big (X_1, Y_1^{\delta }\big) \} \ \Big| \ \big\{N_V=1 \big\} \Big\}\Big) \times  e^{-1}.\end{eqnarray*}
where  $\mathbb{P}(\{N_V=1\}) = e^{-1}$  since $N_V$ is a Poisson random variable of parameter 1.
The random variables $\widetilde X_q$ are i.i.d. uniformly distributed in $V$. Thus,  
\begin{eqnarray*}
 &&   \mathbb{P}\Big(  \Big\{\forall \, 1\le q\le \widetilde N_V, \  \widetilde X_q \not\in B\big (X_1, Y_1^{\delta  }\big) \} \ \Big| \ \big\{N_V=1 \big\} \Big\}\Big) \\
 && \ge  \mathbb{E}\Big ( \Big [1-\frac{\ell\big (B(X_1,Y_1^{\delta   })\big )}{2^{-g(V)}}\Big ]^{\widetilde N_V}\Big ).
\end{eqnarray*}
Observe that, since $\delta>1$,  provided that $g(V)$ is large  enough,  conditionally on $\{N_V\geq 1\}$,  $\ell\big(B(X_1,Y_1^{\delta })\big )\le 2^{-g(V)\delta }$.  This implies that
\begin{eqnarray}\label{aaa}
\mathbb{P}(\widetilde {\mathcal{A}}(U,V,\delta))\ge  e^{-1} \times \mathbb{E}\Big (\big [1-2^{-g(V)(\delta-1 ))}\big ]^{\widetilde N_V}\Big ).
\end{eqnarray}

Let us define $\eta_{g(V)}=2^{-g(V)(\delta-1 )}$. Using that $\widetilde N_V$ is a Poisson random variable of parameter $\widetilde l_V$, a  classical calculus shows that (\ref{aaa}) can be rewritten as 
 $$
\mathbb{P}(\widetilde {\mathcal{A}}(U,V,\delta))\ge e^{-1} e^{-\widetilde l_V \cdot \eta_{g(V)}}.
$$
In order to conclude, it suffices to bound from above the product $\widetilde l_V \cdot  \eta_{g(V)}$. This is achieved by recalling the definition  of $h(V)=\big [(g(V)+1) \delta \big ]+4$, which  implies that 
\begin{eqnarray*}
 \widetilde l_V \cdot  \eta_{g(V)} & \leq & 2^{h(V) -g(V))}  2^{-g(V)(\delta-1 )}\le  16 \cdot 2^\delta.
 \end{eqnarray*}
Thus,  $\widetilde l_V \eta_{g(V)}$ is   bounded from above independently of $V$ by a continuous function of $\de$.  As a conclusion, $\mathbb{P}(\widetilde {\mathcal{A}}(U,V,\delta))$, and thus $\mathbb{P}( {\mathcal{A}}(U,V,\delta))$,  is   bounded from below by some quantity $\kappa_1(\de)$ which is strictly positive and continuously dependent on $\de>1$.  Lemma \ref{bbb} is proved.

\subsection{Diophantine approximation by  rational elements in $\R^d$}  

The question of the validity of $\mathcal{C}_1$ and $\mathcal{C}_2$ may be asked for the rational system in higher dimension $\zu^d$
$$\mathcal{R}^d = \big ((p_1/q,\dots,p_d/q), 1/q^{1+1/d}) \big)_{q\geq 1, 0\le p_i\le q-1}.$$
Again, it follows from Dirichlet's argument   that $ \mathcal{L}_1(\mathcal{R}^d) =[0,1]^d$ (see for instance Theorem 200 in \cite{HW}). The irreducible sub-system of $\mathcal{R}$ consists in the elements  of $ \mathcal{R}$ such that $p_i\land q=1$ for some $i$. 

It is known that $\dim ( \mathcal{L}_\de(\mathcal{R}^d) ) = d/\de$, see \cite{JARNIK2,EGG,BOVEYDODSON}. Using this result, the upper bounds in Theorems \ref{main} and \ref{main2} 
can be proved.
Unfortunately we could not demonstrate neither the weak redundancy property nor $\mathcal{C}_2$ for $\mathcal{R}^d$ (or for any reasonable sub-systems of $\mathcal{R}^d$), so we could not obtain the lower bound.



\begin{thebibliography}{9}


\bibitem{arnold} V.I. Arnold, Small denominators and problems of stability in classical and celestial mechanics. Pspekhi Mat. Nauk {\bf 18},  (1963) 91--192.

\bibitem{MARK} J. Barral, S. Jaffard, N. Fournier, S. Seuret, Markov processes with random singularity spectra. Preprint, 2008.  

\bibitem{BSUbiquity} J. Barral, S. Seuret, Heterogeneous ubiquitous systems in $\mathbb{R}^d$ and Hausdorff dimension. Bull. Braz. Math. Soc. {\bf 38}, (2007) 467--515. 


\bibitem{BSUbiquity2} J. Barral, S. Seuret, Ubiquity and large intersections properties under digit frequencies constraints.   Math. Proc. Cambridge Phil. Soc. {\bf 145} (2008), 527--548.


\bibitem{BS1} J. Barral, S. Seuret, Combining multifractal additive and multiplicative chaos.
Commun. Math. Phys. {\bf  257} (2), (2005) 473-497. 

\bibitem{BS2} J. Barral, S. Seuret, The singularity spectrum of L\'evy processes in multifractal time. 
Adv. Math.  {\bf 214} (1), (2007)  437-468.




\bibitem{Beres} V. Beresnevich, H. Dickinson, S. Velani , {Measure theoretic laws for limsup sets}.  Memoirs of the AMS. {\bf 179} (840), (2006) 1--91. 

\bibitem{Beres3} V. Beresnevich, M.M. Dodson, S. Kristensen, J. Levesley, An inhomogeneous wave equation and non-linear Diophantine approximation. Adv. Math {\bf 217}, (2008) 740--760.

\bibitem{Beres2} {V. Beresnevich, S. Velani}, {A Mass Transference Principle and  the Duffin-Schaeffer conjecture for Hausdorff measures}.  Ann. Math. {\bf  164},  (2006) 971--992. 


\bibitem{BESIC} A.S. Besicovitch, {On the sum of digits of real numbers represented in the dyadic system}. Math. Ann {\bf 110}, (1934) 321--330.


\bibitem{BOVEYDODSON} J. D. Bovey, M. M. Dodson, The Hausdorff dimension of systems of linear
forms. Acta Arith. {\bf 45}  (1986), 337--358.

\bibitem{BUGEAUD} Y. Bugeaud,  Approximation by algebraic integers and Hausdorff dimension.
J. London Math. Soc. {\bf 65}, (2002) 547-559.

\bibitem{BUGEAUD2} Y. Bugeaud,  Sets of exact approximation order by rational numbers.
Math. Ann. {\bf 327},  (2003) 171-190.


\bibitem{BUGEAUD3} Y. Bugeaud,  A note on inhomogeneous Diophantine approximation.
Glasgow Math. J. {\bf 45},  (2003) 105-110. 

\bibitem{DODVEL} M.M. Dodson, M.V, Meli\'an, D. Pestana, S.L.
V\'elani,  Patterson measure and Ubiquity.
Ann. Acad. Sci. Fenn. Ser. A I Math. {\bf 20}, (1995) 37--60.

\bibitem{EGG} H.G. Eggleston,  Sets of fractional dimensions which occur in some problems of number theory. Proc. London Math. Soc. {\bf 54}  (1952), 42--93. 

\bibitem {F2} K. J. Falconer,  Techniques in Fractal Geometry. Wiley, New York (1997).

\bibitem{FalcLI} K. J. Falconer, Sets with large intersection properties. J. London Math. Soc., {\bf 49} (1994), 267--280.

\bibitem{Falc} K. J. Falconer, Representation of families of sets by measures, dimension spectra and Diophantine approximation. Math. Proc. Cambridge Philos. Soc. {\bf128}, (2000) 111--121.

\bibitem{guting} R. G\"uting. On Mahler's function $\theta_1$. Michigan Math. J. {\bf 10},  (1963) 161--179.

\bibitem{HW} G.H. Hardy, E.M. Wright, An Introduction to the Theory of Numbers. Fifth edition, Oxford Science publications, 1979.

\bibitem{HILL} R. Hill, S.L. Velani, The Jarnik-Besicovitch theorem for geometrically finite Kleinian groups.
Proc. Lond. Math. Soc. {\bf 77} (3),  (1998) 524-550. 

\bibitem{HILL2} R. Hill, S.L. Velani, Ergodic theory of shrinking targets.  Invent. math. {\bf  119},  (1995) 175-198. 

\bibitem{JAFFARD} S. Jaffard, The spectrum of singularities of Riemann's function. Rev. Mat. Iberoamericana, {\bf 12}, (1996) 441--460. 

\bibitem{JAFFARD2} S. Jaffard, Old friends revisited: The multifractal
  nature of some classical functions. J. Fourier Anal. Appl.  {\bf
  3} (1), (1997) 1--22.
  
\bibitem{JAFFLACUNARY} S. Jaffard, {On lacunary wavelet series}.  
Ann. Appl. Probab.  {\bf 10} (1), (2000) 313--329.
  
\bibitem{JAFFARD1} S. Jaffard, {The multifractal nature of {L}\'evy
processes}. Probab. Theory Relat. Fields, {\bf 114}, (1999) 207--227.


\bibitem{JARNIK} V. Jarnik, Diophantischen {A}pproximationen und
{H}ausdorffsches {M}ass. Mat. Sbornik {\bf 36},  (1929) 371--381.


\bibitem{JARNIK2} V. Jarnik, Uber die simultanen diophantischen Approximationen. Math. Z. 33
(1931), 503--543.
 
 

\bibitem{khin}   A.Y. Khintchine, Einige S\"atze \"uber Kettenbr\"uche, mit Anwendungen auf die Theorie der diohantischen Approximationen. Math. Ann. {\bf 92}, (1924) 115--125.


\bibitem{khin2}  A.Y. Khintchine, Zur metrischen Theorie der diophantischen Approximationen. Math. Z. {\bf 24} (1), (1926) 706--714.

\bibitem{MATTILA} P. Mattila, Geometry of {S}ets and {M}easures in
{E}uclidian {S}paces. Cambridge Studies in Advanced Mathematics,
Cambridge University Press (1995).

\bibitem{Petro1} G. Petronilho, Global solvavility and simultaneous approximable vectors. J. Differential Equations, {\bf 184} (2002), 48--61.

\bibitem{Petro2} G. Petronilho, Global $s$-solvability, global $s$-hypoellipticity and Diophantine approximation. Indag. Math. (N.S.) {\bf 16}, (2005) 67--90. 

\bibitem{Ptas} B.I. Ptashnik, Improper Boundary Problems for Partial Differential Equations. Nauka Dumka, 1984. 


\bibitem{SLATER} N.B. Slater, Gaps and steps for the sequence $n\theta \mod 1$. Proc. Cambridge Philos. Soc. {\bf 63}, (1967) 1115-1123.

\bibitem{SOS} V.T. S\'os, On the theory of Diophantine approximation, I. Acta Math. Hung. {\bf 8}, (1957) 461--472.  

\bibitem{STRAT} B.O. Stratmann, M. Urbanski,  Jarnik and Julia: a Diophantine analysis for parabolic rational maps,  Math. Scan. {\bf 91},  (2002) 27Ð54.

\bibitem{VEL4}  R.C. Vaughan, S.L. Velani, Diophantine approximation on planar curves: the convergence theory,
  Invent. math. {\bf  116}, (2006) 103-124. 



 \end{thebibliography}
\end{document}